\documentclass{article}

\usepackage{arxiv}

\usepackage{lineno,hyperref}
\modulolinenumbers[5]
\usepackage[dvipsnames]{xcolor}
\usepackage{amsmath}
\usepackage{amsfonts}
\usepackage{graphicx}
\usepackage{amsthm}
\usepackage{graphicx,mathtools}
\usepackage{subcaption}
\usepackage{adjustbox}
\usepackage{bbm}
\usepackage[linesnumbered,ruled]{algorithm2e}
\LinesNumbered 
\usepackage{algpseudocode}
\usepackage{multirow}
\newtheorem{prop}{Proposition}
\newtheorem{remark}{Remark}
\newtheorem{lemma}{Lemma}
\newcommand{\Rbb}{\mathbb{R}}
\newcommand{\Uscr}{\mathcal{U}}
\newcommand{\Lscr}{\mathcal{L}}
\newcommand{\Cscr}{\mathcal{C}}
\newcommand{\Rscr}{\mathcal{R}}
\newcommand{\Oscr}{\mathcal{O}}
\newcommand{\epf}{\hfill \rule[-2pt]{6pt}{6pt}}
\DeclareMathOperator{\argmax}{arg\,max}

\usepackage[english]{babel}

\usepackage{natbib}
\title{Two-stage Robust Optimization Approach for Enhanced Community Resilience Under Tornado Hazards}

\author{
 Mehdi Ansari \\
  American Airlines \\ Fort Worth, TX \\
  \texttt{meansar@okstate.edu}
   \And
 Juan S. Borrero \\
  School of Industrial Engineering and Management\\ Oklahoma State University \\
  \texttt{juan.s.borrero@okstate.edu}
  \And
 Andr\'es D. Gonz\'alez \\
  School of Industrial and Systems Engineering \\ University of Oklahoma\\
  \texttt{andres.gonzalez@ou.edu} \\
}

\begin{document}
\maketitle
\begin{abstract}
Catastrophic tornadoes cause severe damage and are a threat to human wellbeing, making it critical to determine mitigation strategies to reduce their impact. One such strategy, following recent research, is to retrofit existing structures. To this end, in this article we propose a model that considers a decision-maker (a government agency or a public-private consortium) who seeks to allocate resources to retrofit and recover wood-frame residential structures, to minimize the population dislocation due to an uncertain tornado. In the first stage the decision-maker selects the retrofitting strategies, and in the second stage the recovery decisions are made after observing the tornado. As tornado paths cannot be forecast reliably, we take a worst-case approach to uncertainty where paths are modeled as arbitrary line segments on the plane. Under the assumption that an area is affected if it is sufficiently close to the tornado path, the problem is framed as a two-stage robust optimization problem with a mixed-integer non-linear uncertainty set. We solve this problem by using a decomposition column-and-constraint generation algorithm that solves a two-level integer problem at each iteration. This problem, in turn, is solved by a decomposition branch-and-cut method that exploits the geometry of the uncertainty set. To illustrate the model’s applicability, we present a case study based on Joplin, Missouri. Our results show that there can be up to 20\% reductions in worst-case population dislocation by investing \$15 million in retrofitting and recovery; that our approach outperforms other retrofitting policies, and that the model is not over-conservative.
\end{abstract}


\section{Introduction}
Tornadoes are vertical columns of rotating air that are spawned from supercell thunderstorms caused by rotating updrafts that form because of the shear in the environmental wind field \citep{supercell}.
Catastrophic tornadoes are common natural disasters that happen in many populated regions around the world and are particularly concerning in dense high-risk urban areas because of their damage potential. On average, more than 1200 tornadoes happen in the US annually which cause around 60 fatalities \citep{NOAA2022}. Tornadoes have caused damage costs that range between \$183 million to \$9.493 billion per year in the US within the past two decades \citep{tornadocost}. The intensity of tornadoes is  measured by the Enhanced Fujita (EF) scale, which ranks a tornado from 0 (weakest) through 5 (most violent) based on an estimation of their average wind speed. The average wind speed, in turn, is estimated by comparing the observed tornado damages with a list of Damage Indicators and Degrees of Damage  \citep{weather2022}.

Fatalities and injuries from tornadoes have decreased in the past few decades thanks to warning systems based on weather forecasts~\citep{standohar2017vertical,koliou2020development}. For example, the Storm Prediction Center at the National Weather Service keeps a Day 1-8 Convective Outlook for all the US and issues watches or warnings based on the possibility or observation of tornadoes \citep{SPC22}. Tornado warnings are issued as soon as a rotating supercell has been identified by a weather radar; these warnings instruct people in the affected area to take shelter immediately. Even though the warning systems have caused a great reduction in fatalities~\citep{FEMA}, their short lead times (of the order of minutes) do not allow to prevent physical structures from possible destruction. This is further aggravated by the fact that more than 80\% of building stocks in the US are wood-frame buildings that are highly vulnerable to wind damage \citep{van2009performance}. 

Besides weather forecasting,  other alternatives   can be employed to reduce the impact of tornadoes. In the past decade, several studies have shown that retrofitting strategies with simple and inexpensive actions, e.g., improving the roof cover or enhancing the roof sheathing nailing pattern of a house, can improve building codes to make wood-frame buildings more resistant to tornado damage, particularly from damaging tornadoes of EF2 intensity or less~\citep{simmons2015tornado,ripberger2018tornado, masoomi2018wind, koliou2020development, wang2021effect}. In particular, \cite{ripberger2018tornado} shows that a 30\% or more reduction in lifetime damage can be expected by enhancing existing wood-frame buildings with simple retrofitting actions. {Similarly, recent research has studied different recovery strategies and has shown that they have the potential to significantly improve the restoration time and reduce population dislocation~\citep{masoomi2018restoration, farokhnia2020selection, koliou2020development}.




At the federal level in the US, the use of retrofitting has been  identified as a valuable strategy to reduce the impact of natural disasters~\citep{FEMA2021,FEMA2022,WH2022}. At the local level, however, such strategies have not been implemented yet.  In fact, there are federal programs that provide general  guidelines for local administrators (e.g., the emergency management agencies of cities and towns) to design and fund retrofitting plans~\citep{FEMA2021}.  Motivated by this need, the main goal of this paper is to present a model to guide retrofitting and recovery investments in order to enhance the community resilience of an urban area subject to tornado hazards. Here, community resilience is understood as  ``the ability to prepare for anticipated hazards, adapt to changing conditions, and withstand and recover rapidly from disruptions'' \citep{mcallister2015community}. 


Specifically, we assume that there is a decision-maker (a government agency or a public-private consortium; for instance the emergency management department of a city and home insurance companies) that can invest a limited budget to retrofit and recover structures across a geographical area of interest, in order to maximize the community resilience 
due to an uncertain tornado.  We model this decision problem as a two-stage robust optimization model. The first-stage decisions are made before the realization of a tornado and determine what structures have to be retrofitted and to what extent. The second-stage decisions determine recovery strategies for the locations that are damaged by the tornado. Without loss of generality, in this article we consider total \emph{population dislocation}, which is the involuntary movement of people from their residential sites after a severe tornado, as a metric for community resilience to be minimized (therefore, maximizing community resilience would be equivalent to minimizing population dislocation). Our proposed formulation, however, allows for other metrics that assess community resilience.

\looseness-1We formulate the first- and second- stage problems as integer programming problems (IP) whose decision variables select retrofitting strategies and recovery strategies, respectively, for each location.  Both stages share the same budget, which models that, at least a priori, the decision-maker does not know what proportion of the budget should be spent on retrofitting and what proportion on recovery. On the other hand, as an accurate prediction of the location/time where a tornado forms, the motion direction, and its magnitude, is out of reach of the current technology~\citep{TornadoFAQ}, we approach the uncertainty using a robust lens. That is, we assume that given any retrofitting plan, the corresponding worst-case tornado damage happens. In order to have an adequate balance between accuracy and tractability, we model tornado paths as arbitrary line segments on the plane and assume that an area is affected by the tornado if it is sufficiently close to the line segment. This modeling approach results in a mixed-integer non-linear uncertainty set. 

\looseness-1In order to solve the problem, we employ an exact column-and-constraint generation (C\&CG) algorithm based on the method discussed by~\cite{zeng2013solving}. Here, at each iteration, a master problem solves a relaxation of the original problem over a subset of possible tornado scenarios. New scenarios are iteratively generated in a subproblem and added to the master until the objective value does not change. In contrast with previous models that employ the C\&CG algorithm, for instance~~\cite{zeng2013solving,yuan2016robust,velasquez2020prepositioning}, our subproblem is a challenging max-min non-linear integer problem that cannot be solved by standard `dualize and combine' techniques. Thus, the subproblem is solved by a decomposition branch-and-cut algorithm that exploits the geometric properties of the uncertainty set. Particularly, we derive `double' and `triple' conflict constraints to define the master relaxation and we check the feasibility of master solutions by using a \emph{stabbing line algorithm} and a continuous non-convex optimization problem.



Using geographical and population data from~\cite{INCORE},  we use the proposed model to determine optimal retrofitting and recovery actions in Joplin, MO. The results of the case study show that retrofitting actions should be performed in most locations in the geographical center of the city, which are also the locations with higher population density. They also show that spending budget in retrofitting is prioritized over recovery, and that only when a `critical' set of locations in the center of the city are retrofitted, there should be expenditures in recovery.  By performing simulations we also show that the proposed model  outperforms other retrofitting policies, and that the model is not over-conservative: the optimal worst-case population dislocation can be within 10\%  of the maximum population dislocation observed in the simulations.

To summarize, in this paper we make the following contributions:
\begin{itemize}
    \item We propose a novel two-stage optimization model under uncertainty to aid decision-makers in the allocation of resources to retrofit and recover residential wood-frame buildings in tornado-prone regions.
    \item We explicitly model tornado paths as arbitrary line segments on the plane and formulate the problem as a two-stage robust optimization problem with integrality requirements in the first and second stages and in the uncertainty set.
    \item We develop an exact algorithm that embeds a decomposition branch-and-cut algorithm within a column-and-constraint generation method.  By exploiting the geometric properties of line segments, we develop initialization and separation procedures to effectively implement the embedded decomposition branch-and-cut algorithm.
    \item Using real data we provide optimal retrofitting and recovery strategies for Joplin, MO. The results show that there can be up to 20\% reductions in worst-case population dislocation by investing \$15 million; that our approach outperforms other retrofitting policies, and that the model does not suffer from over-conservativeness.
\end{itemize}

The remainder of this paper is organized as follows. Section~\ref{sec:literature review} reviews the literature. Section~\ref{sec:formulation} describes the two-stage robust optimization problem. Section~\ref{sec:solution approach} presents a customized version of the C\&CG algorithm. We describe the decomposition branch-and-cut method to solve the subproblem in C\&CG algorithm in Section~\ref{sec:valid inequality}. In Section~\ref{sec:numerical_expr}, we conduct the numerical experiments on Joplin data as our case study. Lastly, Section~\ref{sec:conclusion} concludes the paper.



\section{Literature review}
\label{sec:literature review}

\subsection{Optimization models to retrofit and recover under tornado hazards}
\looseness-1 The determination of retrofitting and recovery strategies for tornado hazards, particularly using optimization, is a relatively new topic in the literature. \cite{wen2021development}  proposes a multi-objective optimization model to retrofit against tornado hazards. Their model, however, is single-stage,  does not considers recovery,  does not consider uncertain tornado paths, and aggregates all of the uncertainty of the problem into the parameters of the deterministic model. Simulation, on the other hand, is more commonly used to study the behavior of tornadoes, particularly to evaluate the impact of retrofitting strategies, see for example~\cite{strader2016monte, wang2017experimental,masoomi2018restoration,fan2019stochastic, wang2021effect,stoner2021tornado}. These works, however, do not explicitly consider the decision problem of allocating resources to retrofit and recover.

The closest models in the literature to the present work are two-stage robust optimization models for (general) disaster planning, see~\cite{yuan2016robust,ma2018resilience,matthews2019designing,velasquez2020prepositioning,cheng2021robust}. These models, however, are geared towards supply and distribution problems in networks and do not capture the specific patters of the disasters (such as, e.g., paths or ripple effects). Consequently, these models cannot be adapted to deal with tornadoes and the specific retrofitting and recovery setting we study here.  

\subsection{Robust optimization}

Robust optimization methods have been extensively developed in the past two decades to mathematically model and study different systems under uncertainty. Standard robust optimization models assume that the uncertain parameters belong to a convex set \citep{Ben-TalNemirovski99,BertsimasSim04}, which have many applications such as scheduling \citep{lin2004new}, supply chain \citep{ben2011robust}, transportation \citep{yao2009evacuation}, power system problems \citep{xiong2017distributionally}, among others.
Mixed-integer version of uncertainty sets have also been studied. Exact methods to solve a general robust optimization problem with a mixed-integer uncertainty set are similar to the integer L-shaped method \citep{LaporteLouveaux93} in which a master problem iteratively solves a restricted version of the original problem under a few uncertain scenarios and for the solution vector. Here, a subproblem generates a new scenario by finding the worst-case of objective within a mixed-integer uncertainty set \citep{mutapcic2009oracles,ben2015oracle,ho2018oracles,borrero2021modeling,ansari2022robust}.

\cite{ben2004adjustable} extend the scope of standard single-stage robust optimization problems by introducing the adjustable robust optimization methodology, where a part of decision variables must be determined before the realization of the uncertainty set. The rest of the decision variables are chosen when the uncertain parameters are revealed. The two-stage approach provides a less conservative framework to deal with the uncertainty and it can model a broad range of applications in different areas such as transportation \citep{gabrel2014robust}, networks \citep{atamturk2007two}, investment \citep{takeda2008adjustable}, power systems problems \citep{shams2021adjustable,li2021improved}, among others. 

\looseness-1In general, two-stage robust optimization problems are computationally expensive. To efficiently solve them, decomposition methods are widely employed, under the assumption that the second-stage problem is a linear programming problem \citep{thiele2009robust, zhao2012robust, jiang2012benders, gabrel2014robust}. \cite{zeng2013solving} presents the column and constraint generation method (C\&CG) as an alternative to address this class of robust problems. They showed the C\&CG can outperform previous decomposition for many problem settings. The C\&CG approach has become a popular tool to solve two-stage robust optimization in the past decade~\citep{an2014reliable,6812211,6824272, neyshabouri2017two, 7592421, yuan2016robust,matthews2019designing,velasquez2020prepositioning,cheng2021robust}. An important feature of existing applications of the C\&CG method is that the second-stage optimization problem is convex. In contrast, in our problem the second-stage problem is an IP problem and the uncertainty set is mixed-integer non-linear. Therefore, the standard `dualize and combine' approach that is used in the literature to solve the subproblem in the C\&CG does not apply to our case.


\section{Model formulation}
Next, we define the two-stage robust optimization model to minimize the total population dislocation under an uncertain tornado and present a mathematical formulation of the problem.

\label{sec:formulation}
\subsection{Two-stage robust optimization formulation}
Recall that the first stage decisions seek to determine what locations to retrofit before the realization of a tornado. The second stage decisions, that happen after the uncertainty is revealed, select what recovery strategies should be implemented in the locations that are affected by the tornado. 

Formally, suppose $S$ denotes the set of retrofitting strategies and $P$ denotes the set of recovery plans. We consider a set of locations of interest $L$, hereafter referred to as locations for simplicity, which are placed in a 2-dimensional plane. Let $w_{\ell s}$ be the population dislocation estimation pre-tornado at location $\ell \in L$ under retrofitting strategy $s \in S$ (in most cases $w_{\ell s}=0$) and let $g_{\ell sp}$ be the population dislocation post-tornado at location $\ell\in L$ assuming that the retrofitting strategy $s \in S$ is used and that the recovery plan is $p \in P$. 


We assume that the decision-maker has a limited budget of $A\ge 0$ to invest in the retrofitting and recovery plans.
\textcolor{black}{This assumption allows the model to allocate the budget endogenously, rather than having the decision-maker setting the proportions before hand. If required, the model is able to handle fixed budgets by adding a linear constraint in the first-stage problem.} Let $d_{\ell s}$ be the retrofitting cost of all the buildings in location $\ell \in L$ under strategy $s \in S$, and let $c_{\ell sp}$ denote the cost of using recovery plan $p\in P$ on the buildings in location $\ell \in L$ if the retrofitting strategy $s \in S$ is applied. 
\textcolor{black}{All parameters in both stages are computed in Section~\ref{sec:DefParam} in a real example, by using the publicly available data and engineering models for wood-frame buildings in the literature.} 

Define the decision variables $f_{\ell s}$ and $r_{\ell sp}$ as
\[
f_{\ell s} = \begin{cases} 1,  &\text{ if the buildings in location $\ell \in L$ are retrofitted using strategy $s \in S$} \\ 0,  &\text{ otherwise} , \end{cases}
\]
\[
r_{\ell sp} = \begin{cases} 1, &\text{ if the buildings in $\ell \in L$ are retrofitted with strategy $s\in S$}\\
 &\text{ and recovered with plan $p\in P$} \\ 
0, &\text{ otherwise} , \end{cases}
\]
and let $z_{\ell}$, $\ell\in L$, be the binary variables that represent the coverage of a tornado across the locations, that is:
\[
z_{\ell} = \begin{cases} 1 \text{ if location $\ell \in L$ is affected by the tornado} \\ 0 \text{ otherwise} . \end{cases}
\]
Then, the two-stage robust optimization model is defined as follows
\begin{subequations}
\label{eq:two_stage_model}
\begin{align}
    v = \min & \sum_{\ell \in L}\sum_{s\in S} w_{\ell s}f_{\ell s} + \max_{z\in\mathcal{U}} Q(z,f) \label{eq:ts-objfun}     
    \\ \text{s.t. } &
    \sum_{s \in S} f_{\ell s} = 1 & \forall \ell \in L \label{eq:ts-constr1} 
    \\
    & f\in\{0,1\}^{|L||S|},
\end{align}
\end{subequations}
where the second stage problem $Q(z,f)$ for given $z = (z_\ell :\ell \in L)$ and $f=(f_{\ell s}:\ell\in L,s \in S)$ is 
\begin{subequations}
\label{eq:second_stage}
\begin{align}
    Q(z,f) = \min & \sum_{\ell\in L}z_\ell \sum_{s\in S} \sum_{p\in P}g_{\ell sp} r_{\ell sp} \label{eq:ss-objfun}
    \\
    \text{s.t. }& \sum_{\ell \in L} \sum_{s \in S}\sum_{p \in P}c_{\ell sp}r_{\ell sp} \le A - \sum_{\ell \in L}\sum_{s \in S} d_{\ell s}f_{\ell s} \label{eq:ss-constr1}
    \\
    & \sum_{p \in P}  r_{\ell sp} = f_{\ell s} & \forall \ell \in L, s \in S \label{eq:ss-constr2c}
    \\
    &r\in\{0,1\}^{|L||S||P|}\label{eq:ss-constr3}.
\end{align}
\end{subequations}


The first term of the objective function \eqref{eq:ts-objfun} measures the pre-tornado dislocation across all locations due to implementing the retrofitting strategies. The second term measures the population dislocation after the worst-case tornado with respect to the retrofitting strategies selected at the first stage. Constraint \eqref{eq:ts-constr1} ensures that each location picks exactly one retrofitting strategy. We make the assumption that there is a ``do-nothing'' strategy in $S$ with zero cost.

\looseness-1The second stage optimization problem~\eqref{eq:second_stage} selects the recovery plans that minimize the population dislocation given a retrofitting strategy vector $f$ and the tornado coverage vector $z$. The budget constraint \eqref{eq:ss-constr1} ensures that the overall cost for the execution of recovery and retrofitting strategies is $A$. Constraint \eqref{eq:ss-constr2c} forces each location to pick exactly one recovery strategy corresponding to the selected retrofitting strategy. As before, we assume there exist a ``do-nothing'' strategy in $P$ with zero cost. Next, we model tornado paths by ensuring that the decision vector $z$ belongs to the uncertainty set $\mathcal{U} \subseteq \{0,1\}^{|L|}$ which contains all feasible tornado damages over locations in $L$.

\subsection{Uncertainty set formulation}\label{subsec:uncertainy_formul.}

We make three assumptions to model a tornado and its associated damage:
\begin{itemize}
    \item[A.1.] Tornado paths are represented by line segments of at most a  length $E$ in the (real) plane; $E>0$ is a parameter set by the decision maker\textcolor{black}{, and can be easily modified to study the effects of different tornado lengths}.
    \item[A.2.] Each location $\ell\in L$  represents a cluster of buildings that are geographically close. \textcolor{black}{Such clusters can depict individual buildings, blocks, neighborhoods, counties, or any other arrangement of interest to the decision-maker}.
    \item[A.3.] A location $\ell\in L$ is damaged by the tornado if an aggregate measure of the coordinates of all buildings in $\ell$ (e.g., the centroid of all the coordinates of the buildings in that location), is within $\Delta$ units of the tornado path, $\Delta>0$ being a parameter set by the decision-maker.
\end{itemize}

\looseness-1 \textcolor{black}{Assumption A.1. follows from considering historical records of tornadoes, such as evidenced in Figure~\ref{fig:okcHistoricalMap}, which depicts the tornado paths in Oklahoma County between 1950 and 2020. Observe that, with a few exceptions, the paths can be adequately approximated by line segments. This behavior is also observed across the US, not only Oklahoma, see the damage assessment tool in~\cite{damAssesment}. Whereas A.1. might be seen as restrictive, it is a reasonable and uncluttered mathematical representation of a complex weather phenomenon.}
\begin{figure}[h]
    \centering
    \captionsetup{justification=centering}
    \includegraphics[width=6.8cm]{ 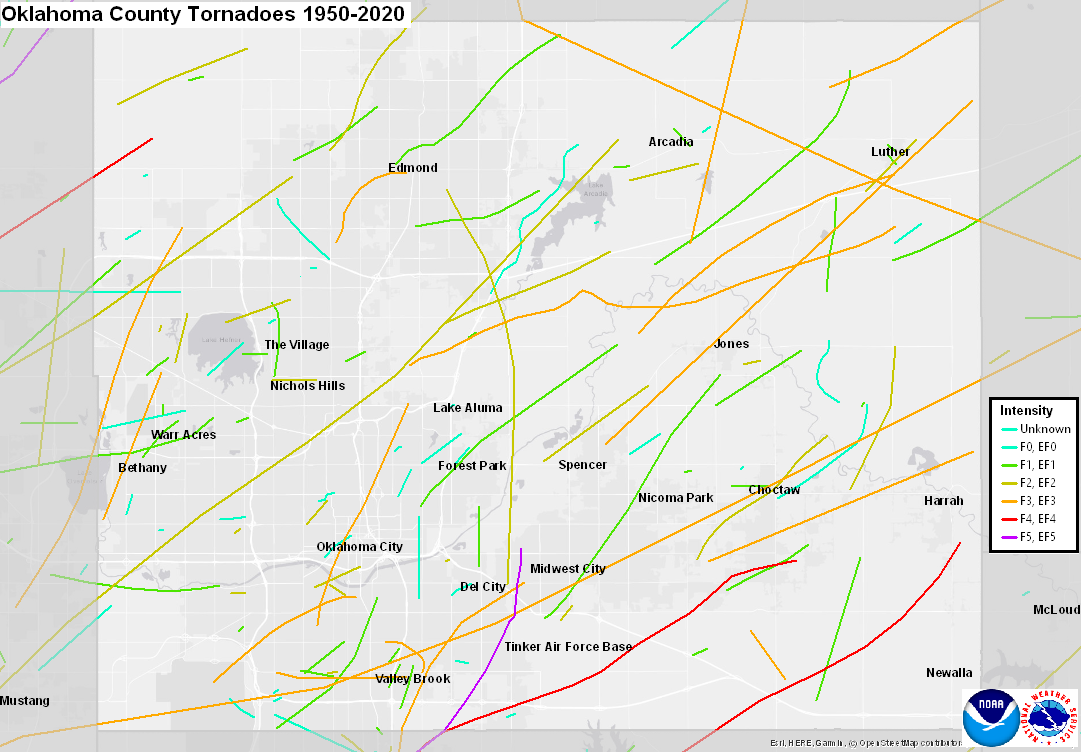}
    \caption{Map of Oklahoma County tornado paths between 1950-2020 \citep{okHist}.}
    \label{fig:okcHistoricalMap}
\end{figure}

\looseness-1\textcolor{black}{Assumption A.2. is made to avoid politically inviable policies. Allowing full granularity (by associating a location with each building) might result in policies that prescribe retrofitting and recovery actions for a given building, but no retrofitting nor recovery to its immediate neighbors. Clearly, such a  policy might raise fairness concerns, particularly if the decision-maker is a public entity. We note that our model works for any clustering approach; all clustering information required by the model is summarized via the parameters of the model.}

\textcolor{black}{Assumption A.3. defines the damage due to the tornado. Considering the variable nature of tornado damage,  $\Delta$ can be seen as an upper-bound on the width observed in historical data for a given area to study the worst-case scenario. Due to the uncertain nature of tornado behavior, defining specific rules for modeling the coverage during a tornado remains a challenging task}.

We assume that the length of tornado $q$ is at most a given parameter $E$, and that the two endpoints of $q$ are located within a sufficiently large rectangle $R = [R_L^{(1)}, R_U^{(1)}]\times[ R_L^{(2)},R_U^{(2)}]\subseteq\Rbb^2$ that contains all locations in $L$. Let $e_0 = (e_0^x,e_0^y) \in R$ and $e_1 = (e_1^x,e_1^y) \in R$ be the endpoints of line segment $q$ (see Figure~\ref{fig:Capsulshape of tornado}). The uncertainty set $\mathcal{U}$ can be defined by the following mixed-integer non-linear programming (MINLP) formulation
\begin{subequations}
\label{eq:uncertain_LS}
    \begin{align}
    \mathcal{U} = \bigg\{z\in\{0,1\}&^{|L|}: \exists e_0,e_1\in R, t \in \mathbb{R}^{|L|},
    \\
    &\label{eq:uncertain_LS_coverage} \|e_0 + t_\ell(e_1 - e_0) - (x_\ell,y_\ell)\|\le \Delta + M_\ell(1-z_\ell), \ \forall \ell \in L,
    \\
    &\label{eq:uncertain_LS_tRange} 0\le t_\ell \le 1, \ \forall \ell \in L,
    \\
    &\label{eq:uncertain_LS_length}\|e_0-e_1\| \le E  \bigg\},
    \end{align}
\end{subequations}
where $\|.\|$ calculates the Euclidean distance between the vector coordinates and where $(x_\ell,y_\ell)$ denote the coordinates of location $\ell\in L$ on the plane. 
\begin{figure}[htb!]
    \centering
    \captionsetup{justification=centering}
    \includegraphics[width=6cm]{ 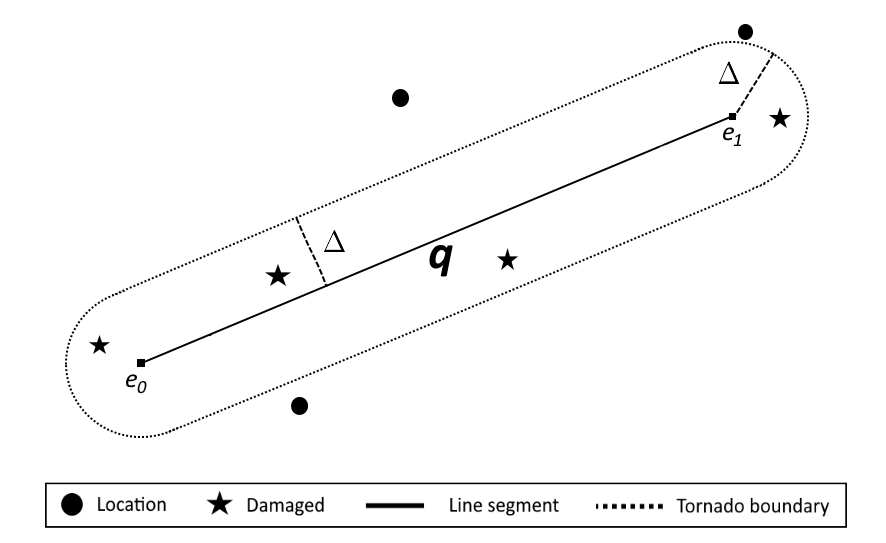}
    \caption{The line segment $q$ represents a tornado central line. Locations within $\Delta$ distance of the line segment (stars) are covered by the tornado.}
    \label{fig:Capsulshape of tornado}
\end{figure}

The value of $\|e_0 + t_\ell(e_1 - e_0) - (x_\ell,y_\ell)\|$ in Constraint \eqref{eq:uncertain_LS_coverage} calculates the distance between a point $e_0 + t_\ell(e_1 - e_0)$ on the line segment $q$ and location $\ell$. This constraint defines $z$-variables such that if $z_\ell=1$ for some $\ell \in L$, then the location $\ell$ is covered by tornado. The value of $M_\ell$ is sufficiently large to ensure the constraint is trivially satisfied if $z_\ell=0$ (the value of $M_\ell$ is set to be the largest value among the Euclidean distances from location $\ell \in L$ to the corners of rectangle $R$). 
Constraint \eqref{eq:uncertain_LS_tRange} enforces the range of continuous variable $t_\ell$ to be between 0 and 1, which indicates that $q$ is a finite segment and not an infinite line. Constraint \eqref{eq:uncertain_LS_length} also ensures that the length of a tornado, which is the Euclidean distance between the two endpoints, does not exceed $E$. Note that the problem maximizes over the $z$-variables and that all of them have non-negative coefficients in the objective function. This observation and Constraint~\eqref{eq:uncertain_LS_coverage} imply that $z_\ell=1$ if and only if location $\ell\in L$ is covered by the tornado path. 

\begin{remark}
    The parameters $\Delta$, $E$, $R_L^{(i)}$, and $R_U^{(i)}$, $i=1,2$, are selected by the decision-maker. Also, we note that $E=\infty$ is a valid choice for the length of the tornado path. In this case we refer to the path as a line rather than by a segment. This assumption represents situations where the path is sufficiently long to traverse the entire region of interest. In this case, constraints~\eqref{eq:uncertain_LS_tRange} and \eqref{eq:uncertain_LS_length} are not necessary in the formulation of $\Uscr$.
\end{remark}

We close this section by noting that problem~\eqref{eq:two_stage_model} belongs to the class of NP-hard problems, see Proposition~\ref{pr:NP-hardness} in Appendix~\ref{proof of Proposition Np-hard}. The proof follows by a reduction from the well-known \emph{knapsack problem} and applies even if there is only one tornado path.


\section{A Solution Method Based on the C\&CG Framework}
\label{sec:solution approach}
We present a method to solve the two-stage robust optimization problem~\eqref{eq:two_stage_model}, based on the C\&CG algorithm in~\cite{zeng2013solving}. First, we provide a linear one-level reformulation of problem~\eqref{eq:two_stage_model} by means of an epigraphic reformulation of the worst-case second-stage objective. To this end, note that $\Uscr$ is a finite set and write it as $\Uscr=\{z^1,\ldots,z^n\}$, where $n\ge 1$ is the cardinality of $\Uscr$ and suppose $r^i$ is the  vector of second-stage variables corresponding to scenario $z^i$. The one-level linear mixed-integer programming (MIP) reformulation of~\eqref{eq:two_stage_model} is:
\begin{subequations}
\label{eq:TS_expanded}
\begin{align}
    v = \min & \sum_{\ell\in L}\sum_{s\in S} w_{\ell s}f_{\ell s} + \theta
    \\
    \text{s.t. }
    &\sum_{s \in S} f_{\ell s} = 1 & \forall \ell \in L
    \\
    & \theta \ge  \sum_{\ell\in L}z^i_\ell \sum_{s\in S} \sum_{p\in P}g_{\ell sp}r_{\ell sp}^i  & \forall i\in\mathcal{I} \label{eq:TS_expanded_secondStage}
    \\
    & \sum_{\ell \in L} \sum_{s \in S} \sum_{p \in P}c_{\ell sp}r^{i}_{\ell sp} + \sum_{\ell \in L}\sum_{s \in S} d_{\ell s}f_{\ell s} \le A & \forall i\in\mathcal{I} \label{eq:TS_expanded_budget}
    \\
    &\sum_{p \in P} r^i_{\ell s p} = f_{\ell s} & \forall \ell \in L, s\in S, i\in\mathcal{I} \label{eq:TS_expanded_uni}
    \\
    &f\in\{0,1\}^{|L||S|}
    \\
    &r^i\in\{0,1\}^{|L||S||P|} & \forall i\in\mathcal{I},
\end{align}
\end{subequations}
where $\mathcal{I}=[n]:=\{1,\ldots,n\}$. Constraint \eqref{eq:TS_expanded_secondStage} implies that the minimum value of $\theta$ is equal to the maximum population dislocation among all tornadoes in $\mathcal{U}$. Constraint \eqref{eq:TS_expanded_budget} and \eqref{eq:TS_expanded_uni} ensure that the recourse vectors $r^i$ corresponding to each scenario, meet the constraints of the second stage.

\begin{remark}
    \looseness-1 Solving problem~\eqref{eq:TS_expanded} directly is not practical in general, given that the uncertainty set $\Uscr$ might have a large number of scenarios in terms of the number of locations. If $E=\infty$, then there is a polynomial number of scenarios in $\Uscr$ which can be constructed based on the stabbing line algorithm in Section~\ref{sec:stabbing_line} (in this case, nevertheless, using~\eqref{eq:TS_expanded} to directly solve the problem remains impractical because there would be at least a quadratic number of scenarios in terms of the number of locations). If $E<\infty$, it remains an open question to determine whether there are polynomially many scenarios in $\mathcal{U}$. 
\end{remark}

We employ a modified version of the C\&CG to solve formulation \eqref{eq:TS_expanded}. 
The algorithm is set to iteratively solve a relaxation of \eqref{eq:TS_expanded} over a subset $\mathcal{I}^k \subseteq \mathcal{I}$ of scenarios at iteration $k$ to obtain an optimal solution $(\theta^k,f^k)$. If given $f^k$, there exists a tornado coverage $z^{i_k} \in \mathcal{U}$ that produces a dislocation greater than $\theta^k$, then we add scenario $z^{i_k}$ to the master problem at the next iteration. Scenario $i_k$ is found by solving the following subproblem, evaluated at $f^k$:
\begin{align}\label{eq:tornado_bilevel_problem}
    \Phi(f) = \max_{z\in \mathcal{U}} \min\left\{\sum_{\ell\in L}z_\ell \sum_{s\in S} \sum_{p\in P}g_{\ell sp} r_{\ell sp}:r \in \mathcal{R}(f)\right\}, 
\end{align}
where $\mathcal{R}(f)$ is the set of feasible solutions of problem $Q(z,f)$; see constraints \eqref{eq:ss-constr1}--\eqref{eq:ss-constr3}. 
Specifically, if the optimal solution is such that $\Phi(f^k)>\theta^k$, then $z^{i_k}$ is the tornado coverage $z$ that attains $\Phi(f^k)$. Note that when the variables and constraints associated with scenario $i_k$ are added at iteration $k+1$, the solution $(\theta^k, f^k)$ is no longer feasible in the master relaxation.

Let $v^k$ be the value of the master problem at iteration $k$, that is, problem~\eqref{eq:TS_expanded} with $\mathcal{I}$ replaced by $\mathcal{I}^k$, with $\mathcal{I}^k=\{i_0,i_1,\ldots,i_{k-1}\}$ and let $\Phi(f^k)$ be the subproblem for the retrofitting strategy $f^k$ found at iteration $k$ of the algorithm. The C\&CG method is presented in Algorithm \ref{alg:algorithm_seq}.

\begin{algorithm}[H]
\footnotesize{
\label{alg:algorithm_seq}
\KwData{Set $\mathcal{U}, L, S$, and $P$}
\KwResult{$f^{k}$}
Set $k = 0$, $\mathcal{I}^1 = \{{i_0}\} \subseteq \mathcal{U}$, $LB=-\infty$, and $UB=\infty$\;
\While{$UB - LB > 0 $}{
Set $k = k + 1$ \; Solve the master MIP $v^k$ and let $f^{k}$ be the optimal retrofitting strategy.
Update $LB = v^k$\;
Solve the subproblem $\Phi(f^k)$ and let ${i_k}$ be the index of the $z$-optimal solution in $\mathcal{I}$.
Update $UB= \min\{UB,  \sum_{\ell\in L}\sum_{s\in S} w_{\ell s}f^k_{\ell s} + \Phi(f^k)$\} and $\mathcal{I}^{k+1} = \mathcal{I}^{k} \cup \left\{{i_k}\right\}$\label{algstep:CCG_Subproblem}\;}
\caption{C\&CG algorithm to solve $v$}}
\end{algorithm}

Algorithm \ref{alg:algorithm_seq} begins with an arbitrary feasible scenario ${i_0} \in \mathcal{I}$. In each iteration $k$, the relaxed MIP problem $v^k$ is solved. The value of $v^k$ then provides a lower bound for the value of $v$. We note that the values of $v^k$, $k\ge 0$, are non-decreasing in $k$. Then, the subproblem $\Phi(f^k)$ generates a new feasible scenario ${i_k}$ for the next iteration and updates the upper bound. The algorithm terminates when the subproblem does not find a scenario that violates the lower bound value, which implies the upper bound is equal to the lower bound. \cite{zeng2013solving} prove that the C\&CG algorithm converges in finite time.

Note that Step~\ref{algstep:CCG_Subproblem} in Algorithm~\ref{alg:algorithm_seq} is expensive computationally because it requires solving the MINLP problem $\Phi(f)$ in~\eqref{eq:tornado_bilevel_problem}. The approach to solve the subproblem in standard applications of the C\&CG  relies on using duality or the Karush Kuhn Tucker (KKT) conditions on the second-stage problem~\citep{zeng2013solving}. Notice that $\Phi(f)$ is a bilevel problem with integrality requirements at both levels. Therefore, we cannot use either of these approaches. Alternatively, we propose a decomposition branch-and-cut method next.

\section{A Decomposition Branch-and-Cut (DBC) Algorithm to Solve $\Phi(f)$}
\label{sec:valid inequality}

The DBC method is based on a one-level reformulation of $\Phi(f)$. We replace the $\min$ in the objective by its epigraphic reformulation and replace the uncertainty set requirements in terms of linear (although exponentially many in the worst-case) constraints. Next, we discuss how to initialize the constraints in the master problem and how to perform the ``two types'' of separations.


\subsection{Overview of the DBC algorithm}\label{subsec:DBCoverview}

Suppose $\mathcal{C}$ is the set of combinations of locations that cannot be covered by a tornado path \textcolor{black}{, i.e., for each $z\notin \Uscr$, there exists a set of locations $C=\{\ell\in L: z_{\ell}=1\} \in \Cscr$}. A one-level linear reformulation of $\Phi(f)$ is given by:
\begin{subequations}
\label{eq:decomposition_problem}
\begin{align}
    \Phi(f) = \max\ & \eta
    \\
    \text{s.t. } & \eta \le \sum_{\ell\in L}z_\ell \sum_{s\in S} \sum_{p\in P}g_{\ell sp} r_{\ell sp} & \forall r \in \mathcal{R}(f) \label{eq:bound_second_stage}
    \\
    & \sum_{\ell\in C} z_{\ell} \le |C|-1 & \forall C\in \mathcal{C} \label{eq:infeasible_comb_indices}
    \\
    & z \in \{0,1\}^{|L|}.
\end{align}
\end{subequations}
The correctness of the reformulation follows from the maximization sense in~\eqref{eq:decomposition_problem}.
\textcolor{black}{Note that $r$ in~\eqref{eq:bound_second_stage} represents a fixed value (not a variable); and that there is such an $r$ for each element of $\Rscr(f)$.}

Observe that $\mathcal{R}(f)$ and $\Cscr$ are finite sets with potentially exponentially many elements, thus $\Phi(f)$ cannot be solved directly by a commercial solver. The DBC starts with a relaxation of $\Phi(f)$ by replacing $\mathcal{R}(f)$ and $\Cscr$ with subsets $\mathcal{R}^0(f)$ and $\Cscr^0$, respectively. The resulting problem is referred as the \emph{master relaxation}. The standard branch-and-cut (BC) algorithm then starts solving the master and prunes nodes by bound and infeasibility. When at some node $h$ an integral solution $z^h$ is found (i.e., a solution that is feasible in the master), its feasibility with respect to the original formulation in~\eqref{eq:decomposition_problem} must be verified.

To this end, Algorithm~\ref{alg:Separation_Procedure} takes the master feasible solution $(\eta^h,z^h)$ at node $h$ of the BC tree and verifies its feasibility: first, it checks whether $z^h$ satisfies all \emph{$\Cscr$-constraints~\eqref{eq:infeasible_comb_indices}}, i.e., it checks if $z^h\in\Uscr$. If $z^h$ does not pass the first check, then let $C^h=\{\ell\in L\colon z_\ell^h=1\}$ be the \emph{active locations} in $z^h$. Observe that $C^h$ is a combination of locations that cannot be covered by a tornado path, thus a cut~\eqref{eq:infeasible_comb_indices} with $C=C^h$ is added on the fly.
If $z^h$ passes the first check, then the algorithm checks whether $z^h$ satisfies all $\Rscr(f)$-constraints~\eqref{eq:bound_second_stage}.  
If $z^h$ also passes the second check, then $z^h$ is feasible in $\Phi(f)$; if not, then a cut~\eqref{eq:bound_second_stage} with $r=r^h$ is added on the fly, $r^h$ being the optimal solution of the second-stage problem associated with $f$ and $z^h$. The performance of Algorithm~\ref{alg:Separation_Procedure} depends on the quality of $\Cscr^0$ and $\Rscr^0(f)$ and on having fast separation routines. We discuss these issues next.

\begin{algorithm}[H]
\footnotesize{
\label{alg:Separation_Procedure}
\KwData{Master feasible solution $(\eta^h,z^h)$}
\KwResult{Generating cuts on the fly, if $(\eta^h,z^h)$ is infeasible}
\uIf{$z^h \notin \mathcal{U}$\label{algStep: check feasibility}}{
Define $C^h = \{\ell \in L: z_\ell^h=1\}$\;
Add constraint~\eqref{eq:infeasible_comb_indices} with $C=C^h$ on the fly\;
}
\uElseIf{$\eta^h> Q(z^h,f)$}{
Define $r^h=\arg\min\{\sum_{\ell\in L}z^{h}_\ell \sum_{s\in S} \sum_{p\in P}g_{\ell sp}r_{\ell sp}: r\in\mathcal{R}(f)\}$\;
Add constraint~\eqref{eq:bound_second_stage} with $r=r^h$ on the fly\;
}
\Else{
$(\eta^h,z^h)$ is feasible.
}
\caption{Separation Procedure for the DBC}}
\end{algorithm}

\subsection{Definition of $\Cscr^0$}\label{sec:Strong valid inequalities}
The set $\Cscr^0$ is constructed by identifying pairs and triples of locations that cannot be covered by a single tornado path. It is worth mentioning that the valid constraints associated with these infeasible cases can also be used to tighten the relaxation of $\mathcal{U}$ by the elimination of potentially many infeasible fractional solutions for $z$-variables; we evaluate the performance of DBC with $\Uscr$ enhanced by these valid cuts in Appendix~\ref{sec:performance-compare}.

\subsubsection{Cuts from infeasible pairs.}
Because the length of a tornado is at most $E$ and only locations within a distance of $\Delta$ from a tornado line segment are covered, a tornado cannot cover two locations $\ell_1,\ell_2\in L$ if they are at a distance of more than $E+2\Delta$. Consequently we have the following result:

\begin{prop}
\label{pr:inf_pair}
Let $\Omega = \left\{(\ell_1,\ell_2)\in L^2: \|(x_{\ell_1},y_{\ell_1})-(x_{\ell_2},y_{\ell_2})\|>2\Delta+E\right\}$ be the set of infeasible pairs. Then any $z\in\Uscr$ satisfies that
\begin{equation}\label{eq:conflict}
    z_{\ell_1}+z_{\ell_2}\le1\qquad \forall (\ell_1,\ell_2)\in\Omega.
\end{equation}
\end{prop}

 The proof of Proposition~\ref{pr:inf_pair} is given in Appendix \ref{proof of Proposition1}.


\subsubsection{Infeasible triple cuts.}
We also identify infeasible triples of locations that cannot be covered by a tornado path. To introduce the infeasible set of triples, we assume $E=\infty$ which implies a line representation for a tornado path. Observe that if a line does not cover a subset of locations $Z\subseteq L$, then there is no segment of the line to cover the locations in $Z$ either. For any locations $\ell_1$ and $\ell_2$, the analysis below partitions the plane $({\mathbb R}^2)$ in two regions, denoted by $R'_{\ell_1,\ell_2}$ and $R''_{\ell_1,\ell_2}$. Region $R'_{\ell_1,\ell_2}$ is the set of points $\ell\in\Rbb^2$ for which there exist a line that covers $\ell_1,\ell_2$ and $\ell$, and $R''_{\ell_1,\ell_2}=\Rbb^2\setminus R'_{\ell_1,\ell_2}$. For simplicity, in the next discussion we interchangeably refer to a point in ${\mathbb R}^2$ by its name, as $\ell\in {\mathbb R}^2$, or by its two-dimensional coordinates, as in $(x_\ell,y_\ell)\in{\mathbb R}^2$.


Observe that a line $\lambda\subset\Rbb^2$ covers both $\ell_1$ and $\ell_2$ if and only if $\lambda\cap B_1\not=\emptyset$ and $\lambda\cap B_2\not=\emptyset$, where $B_i$ is the ball with center $\ell_i$ and radius $\Delta$, $i=1,2$. Let $\Lscr_{\ell_1,\ell_2}$ be the collection of these lines, and let $P_{\ell_1,\ell_2}$ be the set of points that belong to any such line:
\begin{equation}
    P_{\ell_1,\ell_2}=\{\ell\in\Rbb^2\colon \exists\lambda\in\Lscr_{\ell_1,\ell_2}\text{ s.t. } \ell\in\lambda\}.
\end{equation}
For any set $T\subseteq\Rbb^2$ let $T(\Delta)$ be the set of points that are at a distance at most $\Delta$ from $T$, that is
\begin{equation}
    T(\Delta)=\{(x_\ell,y_\ell)\in\Rbb^2\colon \exists (x_t,y_t)\in T\text{ s.t. } ||(x_\ell,y_\ell)-(x_t,y_t)||\le\Delta\}.
\end{equation}
It is readily seen that $R'_{\ell_1,\ell_2}=P_{\ell_1,\ell_2}(\Delta)$.

Now, the set $P_{\ell_1,\ell_2}$ can be characterized by the four lines $\tau_1$, $\tau_2$, $\sigma_1$ and $\sigma_2$ that are tangent to both $C_1$ and $C_2$, where $C_1$ and $C_2$ are the circles enclosing $B_1$ and $B_2$, see Figure~\ref{fig:tangents}. Specifically, let $\tilde P_{\ell_1,\ell_2}$ be the set enclosed by these tangents, see Figure~\ref{fig:tangentsarea}, we next show that $\tilde P_{\ell_1,\ell_2}=P_{\ell_1,\ell_2}$. 

\begin{figure}[htb!]
\centering
    \begin{subfigure}[t]{0.45\textwidth}
         \centering
         \includegraphics[width=0.7\textwidth]{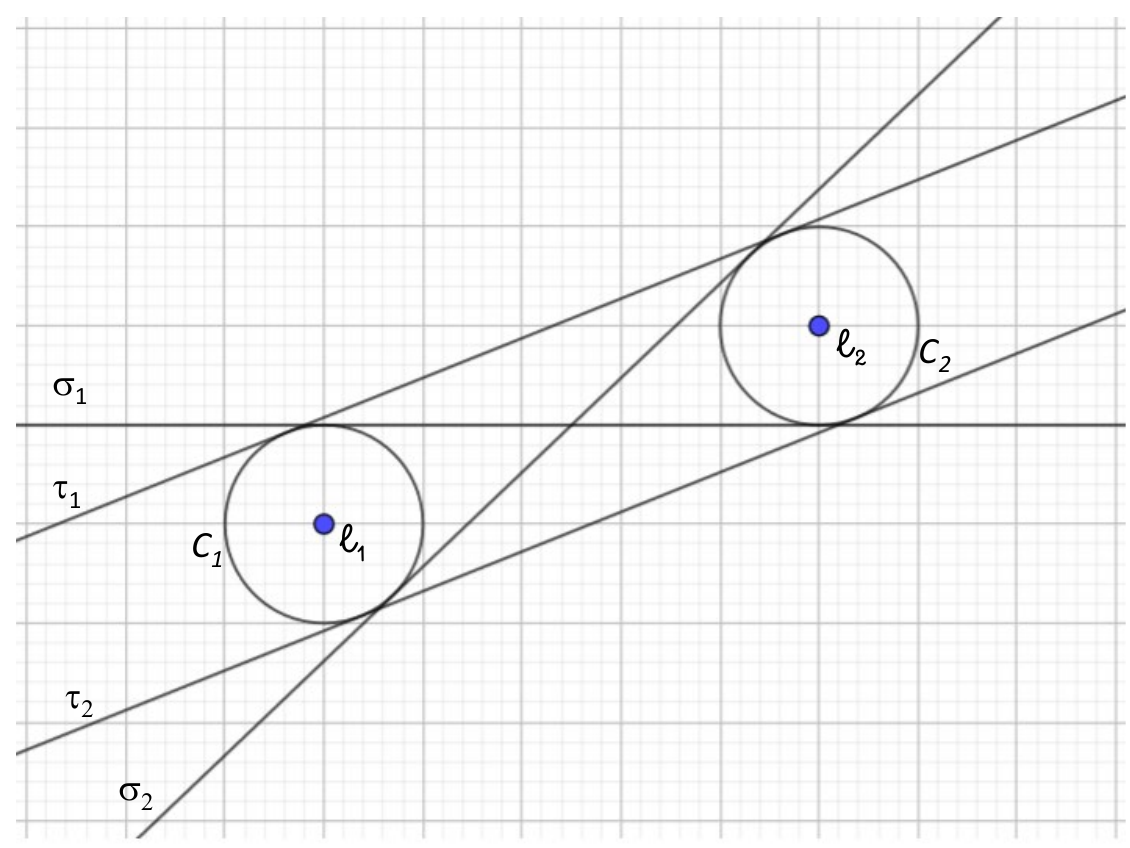}
         \caption{$C_1$ and $C_2$ are circles centered at $\ell_1$ and $\ell_2$, respectively, with radius $\Delta$; $\tau_1$, $\tau_2$, $\sigma_1$, and $\sigma_2$ are the (only) lines that are tangent to both $C_1$ and $C_2$.}
        \label{fig:tangents}
     \end{subfigure}
    \hfill
     \begin{subfigure}[t]{0.45\textwidth}
         \centering
         \includegraphics[width=0.7\textwidth]{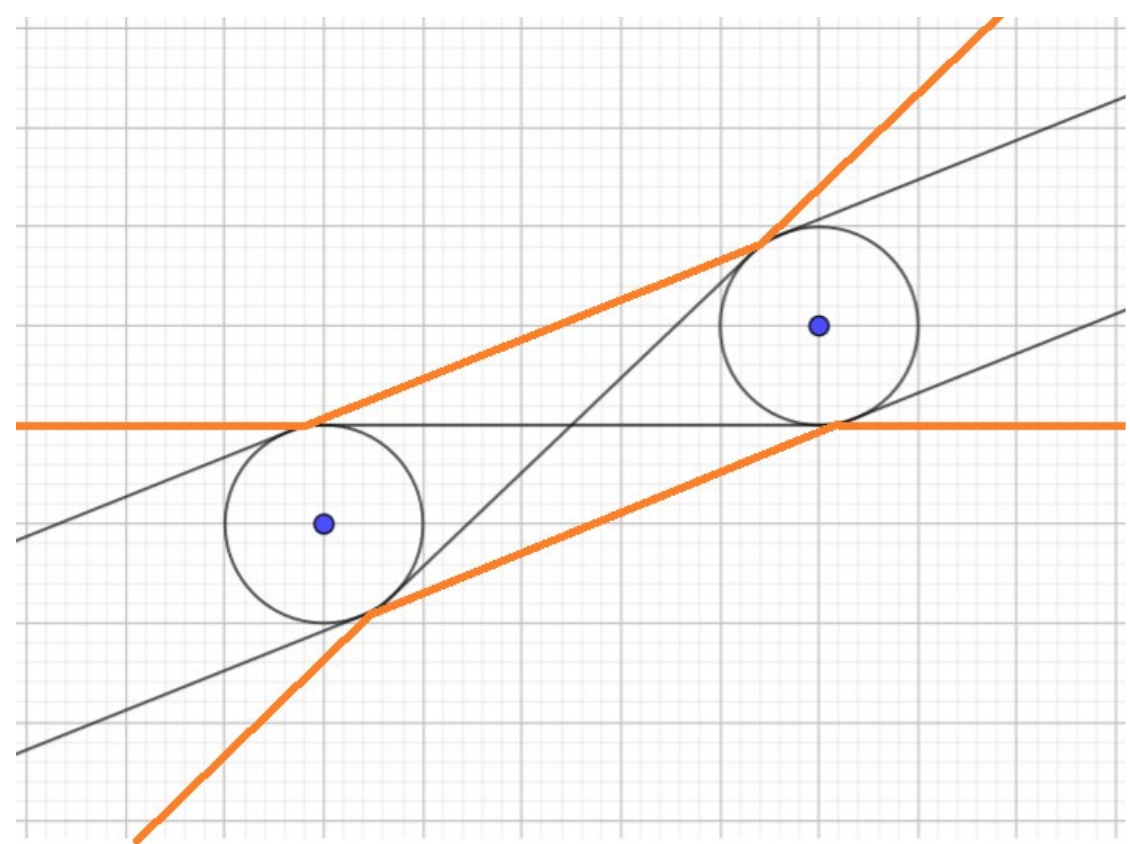}
        \caption{The region $\tilde P_{\ell_1,\ell_2}$ is enclosed by the orange lines.}
        \label{fig:tangentsarea}
     \end{subfigure}
    \caption{{Tangent lines to circles centered at $\ell_1$ and $\ell_2$ with radius $\Delta$}}
    \label{fig:tangentsFigures}
\end{figure}

Indeed, let $\lambda_0$ be the line that passes through $\ell_1$ and $\ell_2$ and suppose that a point $\ell$ is in between $\tau_1$ and $\tau_2$. Then, the line $\lambda$ that passes through $\ell$ and is parallel to $\lambda_0$ intersects both $B_1$ and $B_2$.  Suppose that $\ell$ belongs to the region left of $\ell_1$ between $\tau_2$ and $\sigma_2$. Then, the line $\lambda$ that passes through $\ell$ and the point $\sigma_2\cap C_2$ intersects both $B_1$ and $B_2$. Notice that this argument can be repeated for the remaining regions of $\tilde P_{\ell_1,\ell_2}$, and thus we can conclude that $\tilde P_{\ell_1,\ell_2}\subseteq P_{\ell_1,\ell_2}$. 

Conversely, let $\ell\in P_{\ell_1,\ell_2}$ and let $\lambda\in\Lscr_{\ell_1,\ell_2}$ be the line intersecting $B_1$ and $B_2$ such that $\ell\in\lambda$. If $\lambda$ is one of the tangents, then it is clear that $\lambda\subseteq \tilde P_{\ell_1,\ell_2}$ and it follows that $\ell\in \tilde P_{\ell_1,\ell_2}$. Thus, assume that $\lambda$ is not one of the tangents and let $p_1$ and $p_2$ be two arbitrary elements of $\lambda\cap B_1$ and $\lambda\cap B_2$. Observe that the segment $s$ of $\lambda$ between $p_1$ and $p_2$ must be completely contained in between $\tau_1$ and $\tau_2$. Moreover, as $\lambda$ is not one of the tangents, $s$ intersects both $\sigma_1$ and $\sigma_2$. Because two straight lines can intersect at most once, the previous observations imply that the $\lambda$ cannot intersect $\sigma_1$ nor $\sigma_2$ at any point left of $p_1$ and right of $p_2$; in other words, it must be the case that $\lambda$ is completely contained in $\tilde P_{\ell_1,\ell_2}$, which implies that $\ell\in\tilde P_{\ell_1,\ell_2}$, as desired.

Given the previous considerations, in order to characterize $R'_{\ell_1,\ell_2}$ it is necessary to characterize $\tilde P_{\ell_1,\ell_2}(\Delta)$. From the shape of $\tilde P_{\ell_1,\ell_2}$, it is clear that $\tilde P_{\ell_1,\ell_2}(\Delta)$ looks very similar to $\tilde P_{\ell_1,\ell_2}$, with the difference that its border starts $\Delta$ units up and $\Delta$ units down, see Figure~\ref{fig:region}

\begin{figure}[htb!]
\centering
    \begin{subfigure}[t]{0.45\textwidth}
         \centering
         \includegraphics[width=\textwidth]{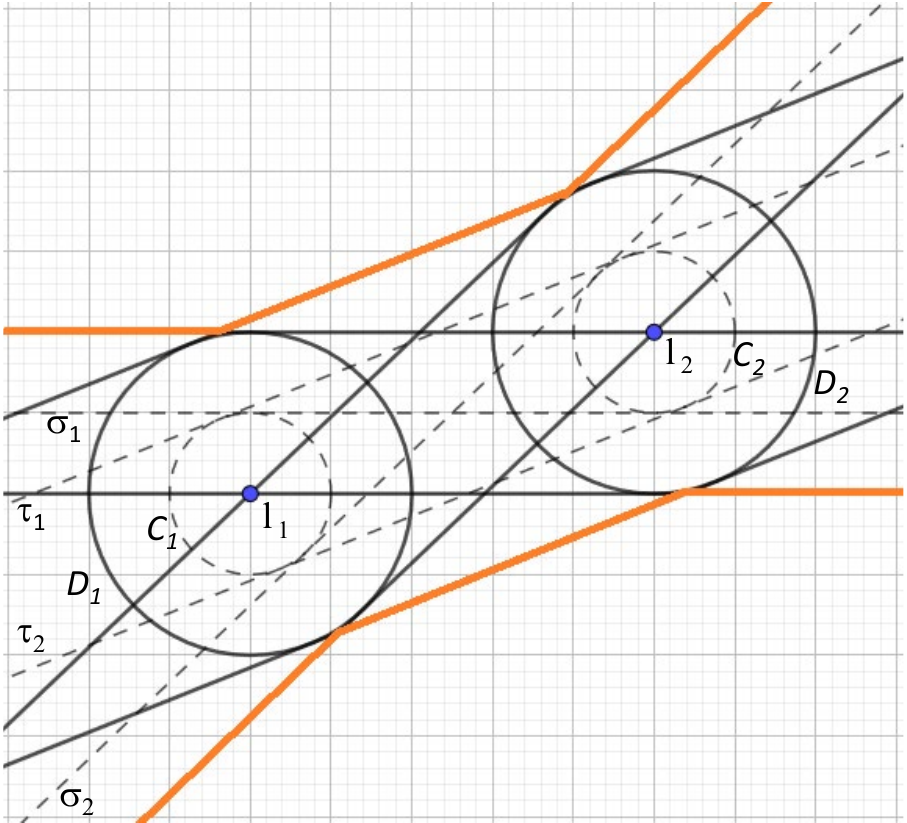}
        \caption{The region $\tilde P_{\ell_1,\ell_2}(\Delta)$ is enclosed by the orange lines. The tangents $\tau_1$, $\tau_2$, $\sigma_1$, and $\sigma_2$, as well as the circles $C_1$ and $C_2$ are shown in dashed lines. The circles $D_1$ and $D_2$ centered at $\ell_1$ and $\ell_2$, respectively, have radius $2\Delta$.}
        \label{fig:region}
        \end{subfigure}
    \hfill
     \begin{subfigure}[t]{0.45\textwidth}
         \centering
         \includegraphics[width=\textwidth]{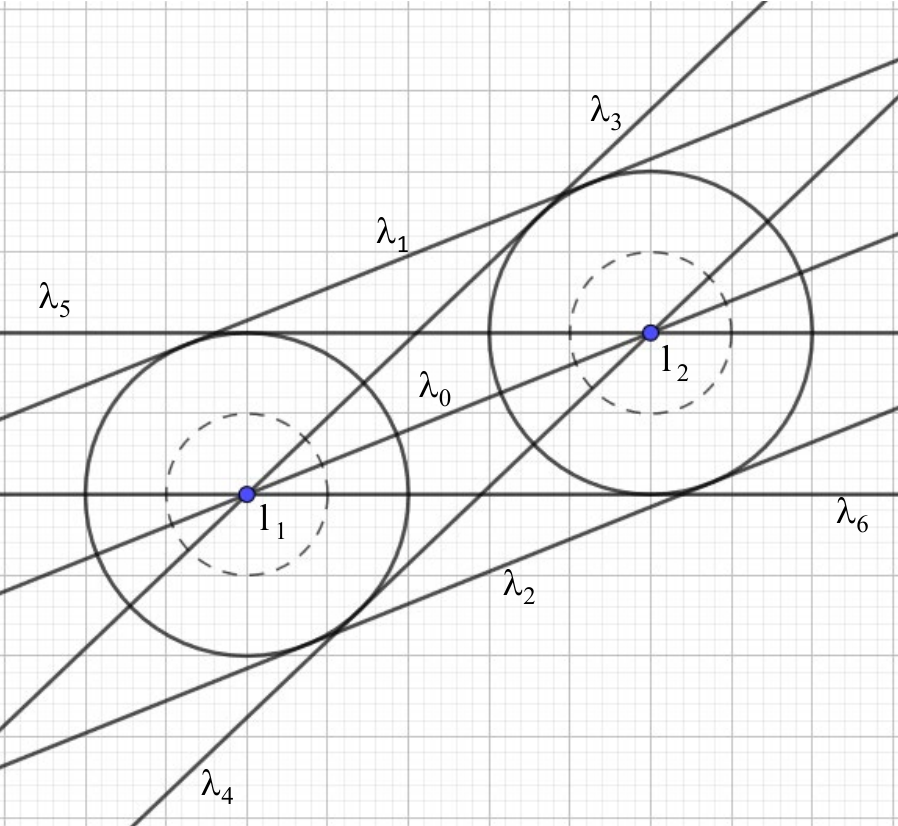}
        \caption{The lines $\lambda_i$, $i=0,\ldots,6$, that define the region $\tilde P_{\ell_1,\ell_2}(\Delta)$. The smaller circles are $C_1$ and $C_2$, the larger circles are $D_1$ and $D_2$.}
        \label{fig:valid_cut}
     \end{subfigure}
    \caption{The region $\tilde P_{\ell_1,\ell_2}(\Delta)$.}
    \label{fig:subspaces-boundaries}
\end{figure}

Next, we use geometric facts to derive an explicit equation for lines defining the boundaries of $\tilde P_{\ell_1,\ell_2}(\Delta)$ for a given $\ell_1$ and $\ell_2$. For ease of exposition hereafter we assume that $x_{\ell_1}\le x_{\ell_2}$ and $y_{\ell_1}\le y_{\ell_2}$. Define $\theta = \arctan{\frac{y_{\ell_2}-y_{\ell_1}}{x_{\ell_2}-x_{\ell_1}}}$ to be the angle of slope of the line $\lambda_0$ which connects $\ell_1$ and $\ell_2$ with respect to the horizontal axis and let $\alpha = \arcsin{\frac{2\Delta}{D(\ell_1,\ell_2)}}$, where $D(\ell_1,\ell_2)$ is the Euclidean distance between $\ell_1$ and $\ell_2$. Then we have that $\widetilde P_{\ell_1,\ell_2}(\Delta)=\{(x,y)\in\Rbb^2\colon$ Eqs~\eqref{eq:P1}--\eqref{eq:P6} hold$\}$, where
\begin{subequations}
\begin{align}
    y&\le \tan{\theta}\cdot (x-x_{\ell_1}) + y_{\ell_1} + \frac{2\Delta}{\cos{\theta}}\label{eq:P1} && (\lambda_1)\\
    y&\ge \tan{\theta}\cdot (x-x_{\ell_1}) + y_{\ell_1} - \frac{2\Delta}{\cos{\theta}}\label{eq:P2}&&(\lambda_2)\\
    y&\le \tan{(\theta+\alpha)}\cdot(x-x_{\ell_1}) + y_{\ell_1}\label{eq:P3}&&(\lambda_3)\\
    y&\ge \tan{(\theta+\alpha)}\cdot(x-x_{\ell_2}) + y_{\ell_2}\label{eq:P4}&&(\lambda_4)\\
    y&\le \tan{(\theta-\alpha)}\cdot(x-x_{\ell_2}) + y_{\ell_2}\label{eq:P5}&&(\lambda_5)\\
    y&\ge \tan{(\theta-\alpha)}\cdot(x-x_{\ell_1}) + y_{\ell_1}\label{eq:P6}&&(\lambda_6).
\end{align}
\end{subequations}
The lines $\lambda_i$, $i=0,\ldots,6$ defined in Equations~\eqref{eq:P1}--~\eqref{eq:P6} are depicted graphically in Figure~\ref{fig:valid_cut}.

Now, let
\begin{equation}\label{eq:defgamma}
    \Gamma_{\ell_1,\ell_2} = \bigl\{\ell \in L\setminus\{\ell_1,\ell_2\}: (x_{\ell}, y_{\ell}) \in \Rbb^2\setminus \tilde P_{\ell_1,\ell_2}(\Delta) \bigr\}.
\end{equation}
The previous discussion is summarized in the following result, which states that if a location belongs to the infeasible set $\Gamma_{\ell_1,\ell_2}$ then at most two locations among $\ell_1,\ell_2$, and $\ell_3$ can be chosen to be covered by a tornado. A more formal proof of this result is given in Appendix~\ref{proof of Proposition2}.

\begin{prop}
\label{pr:triple_validity}
Let $\ell_1$ and $\ell_2$ be two elements of $L$ and suppose that $D(\ell_1,\ell_2)>2\Delta$. Then any $z\in\Uscr$ satisfies that
\begin{equation}\label{eq:triple}
    z_{\ell_1}+z_{\ell_2}+z_{\ell}\le 2\qquad \forall \ell\in\Gamma_{\ell_1,\ell_2}.
\end{equation}
\end{prop}

It is worth mentioning that we do not need to consider all possible combinations of locations to construct the infeasible triple cuts because if $\ell_3 \in \Gamma_{\ell_1,\ell_2}$, then $\ell_1 \in \Gamma_{\ell_2,\ell_3}$ and $\ell_2 \in \Gamma_{\ell_1,\ell_3}$. Also, we note that the pair-wise and triple valid constraints introduced so far do not completely characterize $\Uscr$, as shown in the next remark. However, they can be used to define the initial subset $\Cscr^0 \subseteq \Cscr$ of the master relaxation in the DBC as follows:
\begin{equation}
    \Cscr^0 = \big\{(\ell_1,\ell_2) \in \Omega\big\} \cup \big\{(\ell_1, \ell_2, \ell_3):  \ell_3 \in \Gamma_{\ell_1,\ell_2}, \ell_1,\ell_2\in L, D(\ell_1,\ell_2)>2\Delta\big\}.
\end{equation}

\begin{remark}
\label{remark:counterexample_cuts} 
Define $\mathcal{V}=\{z\in\{0,1\}^{|L|}: \sum_{\ell\in C} z_{\ell} \le |C|-1, \forall C\in \Cscr^0  \}$. The set $\mathcal{V}$ does not necessarily determine all the feasible solutions in the uncertainty set $\mathcal{U}$; specifically $\mathcal{U} \subset \mathcal{V}$ and in general there exist points in $\mathcal{V}$ that do not belong to $\Uscr$. For example, consider a case of three locations $\ell_1,\ell_2$, and $\ell_3$ with coordinates $(0,0),(4,0)$, and $(2,1.1)$ on a plane under a tornado with the maximum length of $E=2$ and the width of $\Delta=1$ (See Figure \ref{fig:counterexample1.jpg}). It can be easily observed that there exists no tornado path that covers all three locations because the only tornado path covering $\ell_1$ and $\ell_2$ is the line segment with endpoints $(1,0)$ and $(3,0)$. Obviously, the location $(2,1.1)$ is not covered by this path. However, $z_{\ell_1}=z_{\ell_2}=z_{\ell_3}=1$ is a feasible solution in $\mathcal{V}$ because it satisfies both the infeasible pair and triple cuts.
\begin{figure}[h!]
    \centering
    \includegraphics[width=5cm]{ 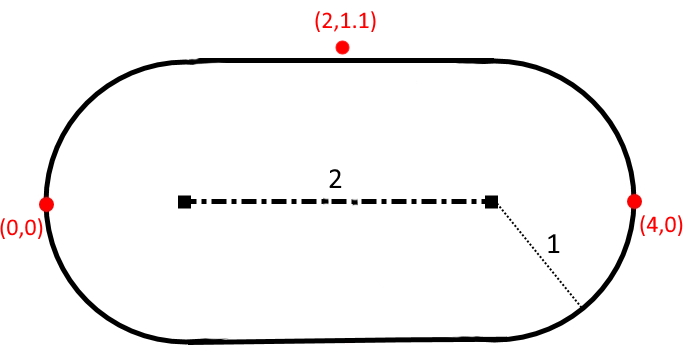}
    \caption{The only possible tornado path covering buildings $(0,0)$ and $(4,0)$.}
    \label{fig:counterexample1.jpg}
\end{figure}
\end{remark}

\subsection{Feasibility check in $\Cscr$}\label{sec:feasibility_check}


Next, we present a method to determine whether a master feasible solution $z^h$ belongs to $\Uscr$. We start our analysis for the case when $E=\infty$ and then present a method for the case when $E<\infty$.

\subsubsection{Feasibility check for $E=\infty$ using the stabbing line algorithm (SLA).}\label{sec:stabbing_line}

\looseness-1The stabbing line problem consists of a collection of (disjoint) circles $\Oscr=\{O_{1}, O_{2}, \ldots, O_{n}\}$ and the objective is to find the line on the plane that intersects the most circles. The SLA is an algorithm that solves this problem and reports $c(\Oscr)$: the maximum number of circles in $\Oscr$ that can be intersected with a line. Consider the collection of circles $\Oscr^h$ with centers at the locations in $C^h$ and radius $\Delta$. Observe that a tornado path with $E=\infty$ covers all locations in $C^h$ if and only if the $c(\Oscr^h)=|C^h|$ (see Figure~\ref{fig:stabbingLine}). Therefore, if $c(\Oscr^h)<|C^h|$, then $z^h\not\in\mathcal{U}$ whereas if $c(\Oscr^h)=|C^h|$ then $z^h\in\Uscr$.

\begin{figure}[h]
    \centering
    \captionsetup{justification=centering}
    \includegraphics[width=6cm]{ 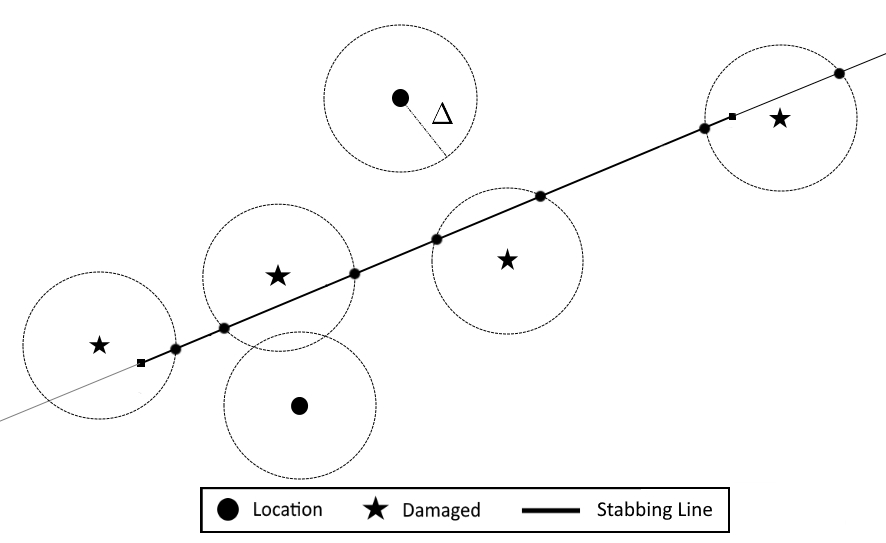}
    \caption{Finding the locations that are covered by a tornado within distance $\Delta$ is equivalent to finding a line that stabs the circles centered at the locations with radius $\Delta$.}
    \label{fig:stabbingLine}
\end{figure}

We explain next the SLA following~\cite{stabline22}. Assume first a set of disjoint circles $\Oscr=\{O_{1}, O_{2}, \ldots, O_{m}\}$, $m>0$. Given any two disjoint circles $O_a$ and $O_b$, find the two external and two internal tangent lines to each pair of circles $(O_a, O_b)$ and record the intersection points on both circles; see Figure~\ref{fig:tangentLines}. Let $p_{ab}^E,q_{ab}^E$ and $p_{ab}^I,q_{ab}^I$ be the intersection points of these external and internal tangent lines, respectively, with $O_a$. Note that a tangent line to $O_a$ that touches any point $p$ in the arc $[p_{ab}^E,p_{ab}^I]$ (or $[q_{ab}^I,q_{ab}^E]$) intersects $O_b$, see Figure~\ref{fig:tangentLines}.
\begin{figure}[htb!]
    \centering
    \captionsetup{justification=centering}
    \includegraphics[scale=0.20]{ 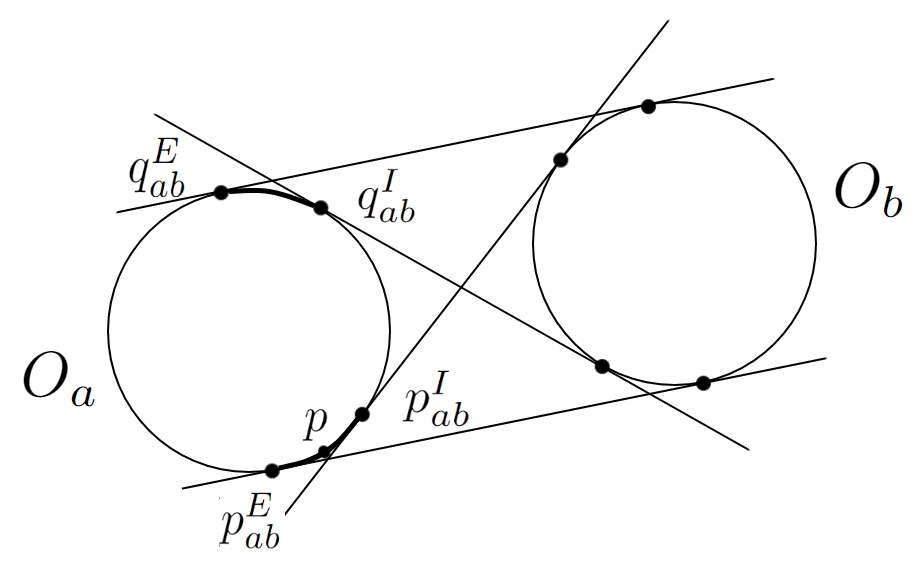}
    \caption{Four tangent lines and the intersection points for a pair of circles $(O_a,O_b)$ are shown. A tangent line at any point $p$ on the arc $[p_{ab}^E,p_{ab}^I]$ also intersects $O_b$.}
    \label{fig:tangentLines}
\end{figure}

Given $O_a$, the SLA computes the arc-intervals $[p_{ab}^E,p_{ab}^I]$ and $[q_{ab}^I,q_{ab}^E]$ for all $b\not=a$; note there exist exactly $2(m-1)$ such intervals for $O_a$. Let $\pi_a$ be an endpoint of any such arc-intervals that belong to the maximum amount of arc-intervals. Then, the tangent line that touches $O_a$ at $\pi_a$ is a line intersecting $O_a$ that intersects the most circles in $\Oscr$. Iterating this process starting from each $a=1,\ldots,m$, will determine a line that intersects the most circles. The SLA can be implemented $O(m^2 \log m)$ time; see \cite{stabline22} for further details on the algorithm.

\looseness-1In our case the circles in $\Oscr^h$ might be non-intersecting, as the locations in $C^h$ might be at a distance less than $2\Delta$. We thus adapt the SLA, noting that if $O_a\cap O_b\not=\emptyset$ then the internal tangents are not necessary and that doing the analysis with the arc-interval $[p_{ab}^E,q_{ab}^E]$ is sufficient for this case.

\subsubsection{Feasibility check if $E<\infty$.}

\looseness-1In this case, if $c(\Oscr^h)<|C^h|$, then we can again conclude that there is no finite line segment that covers $C^h$ and thus $z^h\not\in\Uscr$. Otherwise, if $c(\Oscr^h)=|C^h|$, then there is an infinite line that covers all locations in $C^h$, but there is no guarantee that a line segment of \emph{length at most $E$} covers all locations of $C^h$. Thus, we complete the feasibility check in two steps: first, we use the line resulting from the SLA and projection ideas to determine whether the line can be trimmed to be of length $E$ without compromising its coverage, the details are  in Appendix~\ref{ssec:feasibility}. If even after this step the feasibility of $z^h$ cannot be guaranteed, then we check the feasibility of $z^h$ algebraically using a simplified version of $\Uscr$.

Specifically, consider the set $\tilde\Uscr^h$ defined as
\begin{subequations}
\label{eq:reduced_uncertain_LS}
    \begin{align}
    \tilde\Uscr^h = \bigg\{(e_0,e_1)\in & R^2: \exists t \in [0,1]^{|C^h|}, \|e_0-e_1\| \le E 
    \\
    &\label{eq:reduced_uncertain_LS_coverage} \|e_0 + t_\ell(e_1 - e_0) - (x_\ell,y_\ell)\|\le \Delta, \ \forall \ell \in C^h \bigg\}.
    \end{align}
\end{subequations}
Observe that $z^h\in\Uscr$ if and only if $\tilde \Uscr^h\not=\emptyset$ and  that~\eqref{eq:reduced_uncertain_LS} is far simpler than $\Uscr$ because it has fewer variables and does not involve binary variables nor `big-M' values. Checking whether $\tilde\Uscr^h$ is non-empty can be done using a quadratic (non-convex) continuous optimization solver. The such check can be significantly speed-up by adding linear cuts in $\tilde\Uscr^h$, see the details in Appendix~\ref{ssec:feasibility}.

\subsection{Definition of $\Rscr^0(f)$ and feasibility check for $\mathcal{R}(f)$}
To initialize the master problem we assume that $\Rscr^0(f)=\{r^0\}$ where $r^0$ is an optimal solution for the second-stage problem $Q(z,f)$  assuming no location is hit, i.e., $r^0$ is a solution of $Q(\mathbf{0},f)$. On the other hand, once it has been checked that $z^h\in\Uscr$, we check that $z^h$ satisfies all $\Rscr(f)$-constraints in~\eqref{eq:bound_second_stage} by solving the second-stage problem evaluated at $z^h$, i.e., by solving $Q(z^h,f)$. Let $r^h$ be an optimal solution of $Q(z^h,f)$. If $\eta^h>Q(z^h,f)$ then $(\eta^h,z^h)$ is not feasible in~\eqref{eq:bound_second_stage} and we add the constraint $\eta\le \sum_{\ell\in L}z_\ell\sum_{s\in S}\sum_{p\in P}g_{\ell s p}r_{\ell s p}^h$ at all active nodes of the BC tree. Conversely, if $\eta^h\le Q(z^f,f)$, then $(\eta^h,z^h)$ is a feasible solution of~\eqref{eq:decomposition_problem} and node $h$ is pruned by feasibility after updating the lower bound of the BC. Observe that the initialization and separation procedure involve solving the $NP$-hard problem $Q(z,f)$; however, this problem is relatively simple and can be solved quickly using commercial MIP solvers, see our computational results in Section~\ref{sec:numerical_expr}.

\section{Case study: retrofitting/recovery residential buildings in Joplin, MO}\label{sec:numerical_expr}

Next, we perform several numerical experiments with the proposed model. Our objectives are to study the performance of the two-stage optimization problem, to analyze how the optimal solutions change with the length of the path, to gather high-level insights into the optimal retrofitting and recovery activities, and to compare the optimal policies with other benchmark policies. The experiments are based on location and cost data from Joplin, MO. Joplin was the location of one of the most severe tornadoes in US history that occurred on May 22, 2011. The tornado was rated EF-5 with an average speed of more than 200 mph and crossed 22.1 miles on the ground with a width of one mile. This catastrophe caused the loss of 161 lives and more than a thousand injuries. The total damage was approximated at 3 billion dollars as around 7500 residential structures were destroyed \citep{kuligowski2014final}. 

\subsection{Definition of Parameters}\label{sec:DefParam}
The geographic and demographic data are based on information from Joplin and extracted from the online platform \cite{INCORE}. Particularly, this dataset divides buildings in about 1500 ``block IDs'' (which are contiguous blocks of the city) and includes over 20,000 buildings. For our experiments, we use only the single or multi-story residential wood-frame building archetypes. 

\looseness-1We assume an EF-2 tornado with a wind speed of 134 mph, a width of 0.75 miles (thus $\Delta$=0.75/2), and a length of 5 and $\infty$ miles. Following \cite{wang2021effect}, we consider the retrofitting strategies presented in Table~\ref{tab:retro}; these strategies are analyzed in~\cite{masoomi2018wind} with respect to the component fragility of wood-frame residential building archetypes. \textcolor{black}{The retrofitting activities in Table~~\ref{tab:retro}, have minimal impact in the livability of the buildings; thus we assume that the first stage population dislocation parameter is zero for all locations and strategies}, that is, $w_{\ell s} = 0$ for all $\ell \in L, s \in S$. \textcolor{black}{This assumption, however, does not change the complexity of solving the problem, and if necessary could be relaxed to better depict a particular context or instance of interest}.
\begin{table}[htb!]
\centering
\caption{\label{tab:retro}Retrofitting strategies for the experiments.}
\footnotesize{
\begin{tabular}{c|p{10cm}}
    Strategy $s$ & Description\\
    \hline
    $R_1$ &  Roof covering with asphalt  shingles, using 8d nails to attach roof deck sheathing panel to rafters spaced at 6/12 inch, and the selection of two 16d toenails for roof-to-wall\\
    $R_2$ & Roof covering with asphalt  shingles, using 8d nails to attach roof deck sheathing panel to rafters spaced at 6/6 inch, and the selection of two H2.5 clips for roof-to-wall\\
    $R_3$ & Roof covering with clay tiles, using 8d nails to attach roof deck sheathing panel to rafters spaced at 6/6 inch, and the selection of two H2.5 clips for roof-to-wall.\\
    \hline
\end{tabular}}

\end{table}

The cost of implementing these strategies are computed from the data of ~\cite{INCORE}.
As mentioned in Section~\ref{sec:formulation}, a `do-nothing' strategy with zero cost is also included in the retrofitting strategy set $S$, that is, $S$=\{do-nothing, R1, R2, R3\}. On the other hand, only two recovery strategies are assumed in our experiment, namely, recover ($p=1$) and do-nothing ($p=0$); thus $P=\{0,1\}$. 

The second stage population dislocation parameters, $g_{\ell sp}$, are defined as the expected population dislocation $60$ days after the tornado \textcolor{black}{(Appendix~\ref{sec:30daysRecoveryResults} has additional experiments for 30 days of recovery)}. These values are evaluated based on the \emph{state of damage $d\in D$} immediately after the tornado hits,  where $D$ is the set of damage states. Let $X_\ell$ be the random variable denoting the functionality of $\ell$ after 60 days; where $X_\ell=0$, only if $\ell$ reaches 100\% functionality after $60$ days (no dislocation) and $X_\ell=1$ otherwise (dislocation). Also, let $Y_\ell\in D$ be the random variable representing the damage at $\ell$ immediately after the tornado, and let $\sigma_\ell$ denote the retrofitting strategy. Then, the expected population dislocation after recovery is
\begin{equation}
    g_{\ell s1}=N_{\ell}\times E[X_\ell|\sigma_\ell=s],
\end{equation}
where $N_{\ell}$ is the population of location $\ell \in L$. Since the value of $X_\ell$ depends on $Y_\ell$, using conditional probability and the (reasonable) assumption that $X_\ell$ is independent of $\sigma_\ell$ given the knowledge of $Y_\ell$, we have that
\begin{equation}
    E[X_\ell|\sigma_\ell=s] = \sum_{d\in D} P[X_\ell=1|Y_\ell=d]P[Y_\ell=d|\sigma_\ell=s].
\end{equation}
Following \cite{koliou2020development}, the damage states are given by $D=\{$Minor, Moderate, Extensive, Complete$\}$. We compute $P[X_\ell=u|Y_\ell=d]$ using the repair time distributions described in~\cite{hazus2003multi}. These distributions are lognormal and their parameters can be found in~\cite{koliou2020development}. The probabilities $P[Y_\ell=d|\sigma_\ell=s]$ are computed from the building-level tornado fragility-curves for wood-frame residential buildings in~\cite{masoomi2018wind,koliou2020development}. 

The cost of recovery $c_{\ell s 1}$ is the expected percentage of the total location replacement cost: suppose $R_{\ell}$ is the area of location $\ell \in L$ and $\alpha$ is the cost of replacement per $m^2$. For the wood residential archetype, we set $\alpha=\$862.0 / m^2$ and, the percentage of the location replacement cost to be $r_{minor}=0.5\%$, $r_{moderate}=2.3\%$, $r_{extensive}=11.7\%$, and $r_{complete}=23.4\%$, see \cite{koliou2020development}. The following formula is used to compute the cost of recovery for location $\ell$ if it is retrofitted by strategy $s$:
\begin{equation}
    c_{\ell s 1} = \alpha \times R_{\ell} \times \sum_{d \in D} r_d P[Y=d|\sigma=s].
\end{equation}
We assume the decision-maker takes charge of all recovery expenses. 

On the other hand, if a location does not receive resources from the decision-maker, then the population dislocation will be higher and depends on the resident's actions. We assume the population dislocation for the do-nothing strategy is a factor of the population dislocation resulting from the recovery actions of the decision-maker, i.e.,
\begin{equation}
    g_{\ell s 0} = \mu_{\ell s} g_{\ell s 1} \quad \forall \ell \in L, s \in S.
\end{equation}
The factor $\mu_{\ell s}$ can vary such that $g_{\ell s 0} \in [g_{\ell s 1},N_{\ell}]$: $g_{\ell s 0} = g_{\ell s 1}$ implies that the residents pay a full amount of expenses by themselves to recover their buildings as quickly as the decision-maker, and $g_{\ell s 0} = N_{\ell}$ means that the residents take no recovery action on their own and therefore the whole population is dislocated from location $\ell$. For our experiments, we fix $g_{\ell s 0} = (g_{\ell s 1}+N_{\ell})/2$ which is the midpoint in the range. Note that the cost of the do-nothing strategy is $c_{\ell s 0} = 0, \forall \ell\in L, s \in S$ from the decision-maker's point of view.

\textcolor{black}{We remark that various parameters of our case study we selected rather arbitrarily; for instance, the assumption that recovery means that the building recovers 100\% functionality; that the last day of recovery is 60 days; or that `half' the residents recover their property without assistance. These are reasonable values that are selected for illustration purposes and can be changed by the decision-maker as desired. Generally, modifications of these parameters should not impact the model's performance, and as long as there is data available that account for the modifications, our model should be able to run; see for example Appendix~\ref{sec:30daysRecoveryResults} where experiments with a recovery time of 30 days are reported.}

\subsection{Robust retrofitting and recovery strategies in Joplin, MO}\label{sec:results Robust strategies in Joplin}
We implement the proposed two-stage optimization model to find robust retrofitting and recovery strategies for locations in Joplin. All algorithms are coded in Python 3.9 and solved using Gurobi Optimizer 9.1.1 on a 64-bit Windows operating system with Intel(R) Core(TM) i7-8550U CPU and 8.00 GB of RAM under a time limit of 3600 seconds. The source code is available in online repository \cite{github}. We use Algorithm~\ref{alg:algorithm_seq}, along with the DBC method described in Section~\ref{sec:valid inequality} with the setup that yields the best performance; see our comparative experiments in Table~\ref{Tab1:comparison} in Appendix~\ref{sec:performance-compare}.


In order to define the locations in $L$, we group Joplin ``block IDs'' in 100 center locations using the $k$-means method \citep{macqueen1967classification} and aggregate the data accordingly. Particularly, we assume that when the tornado hits a location then all of its buildings are hit, and that the tornado hits a location only if the location's centroid is within $\Delta$ units of the path. Table~\ref{Tab2:JoplinRobust} compares the population dislocation for when the tornado length is  $E=5$ miles or $\infty$ and the budgets are $A=0,15,30$ million USD. Columns 1 and 2 in Table~\ref{Tab2:JoplinRobust} show $E$ and $A$. Column 3 has the optimal population dislocation under the worst-case tornado for each set of parameters. Columns 4 and 5 report the CPU run time of C\&CG Algorithm~\ref{alg:algorithm_seq} in seconds and its number of iterations, respectively. Column 6 shows the CPU runtime of the DBC method in seconds to solve subproblem $\Phi(f)$. Column 7 compares the CPU time that the DBC algorithm spends in a callback function to verify the feasibility of integral solutions. Column 8 shows the total CPU time to find a feasible solution in the set $\Tilde{\mathcal{U}}^h$. Columns 9 and 10 present the proportion of the total budget that is spent for retrofitting and recovery, respectively.

\begin{table}[!ht]
    \centering
    \captionsetup{justification=centering}
    \caption{Solving the two-stage robust optimization model with different parameters of tornado length and available budget,  \textcolor{black}{after 60 days of recovery}  for 100 locations group in Joplin, MO.}
    \begin{adjustbox}{max width=\textwidth}
    \begin{tabular}{c|c|c|cc|ccc|cc}
    \hline
        Length & Budget & Population  & CCG  & CCG & DBC  & Callback  & $\Tilde{\mathcal{U}}^h$ feasibility & Retrofitting cost & Recovery cost\\ 
        & (\$) & dislocation & run time (sec.) & iteration  & run time (sec.) & run time (sec.) & run time (sec.) & /budget & /budget  \\
        \hline
        \multirow{3}{*}{5} & 0 M & 16318 & 127 & 2 & 126 & 3 & 1 & - & - \\
        ~ & 15 M & 13408 & 307 & 3 & 306 & 11 & 6 & 0.6 & 0.4 \\
        ~ & 30 M & 13091 & 300 & 3 & 300 & 11 & 5 & 0.3 & 0.7 \\ \hline
        \multirow{3}{*}{$\infty$} & 0 M & 17293 & 568 & 2 & 568 & 3 & 0 & - & - \\
        ~ & 15 M & 14236 & 716 & 2 & 716 & 5 & 0 & 0.6 & 0.4 \\
        ~ & 30 M & 13839 & 986 & 3 & 985 & 5 & 0 & 0.4 & 0.6 \\ \hline
    \end{tabular}
    \end{adjustbox}
    
       \label{Tab2:JoplinRobust}
\end{table}

Table~\ref{Tab2:JoplinRobust} shows that the solver is able to find optimal solutions for every experiment in at most 986 seconds. It also shows that solving the subproblems is the most time-consuming step (which is virtually 99\% of the total runtime). In general, the algorithm solves the subproblem for 5 miles line segment faster than the (infinite) line; this fact can be justified by the effect of the conflict constraints~\eqref{eq:conflict} which are not present when $E=\infty$. In addition, for all experiments, the callback function (i.e., Algorithm~\ref{alg:Separation_Procedure}) is very fast, which implies that the bottleneck of the entire algorithm is the branch-and-bound that solves the relaxation of the subproblem. As expected, for the infinite line, the results show zero seconds to check the feasibility using $\tilde{\mathcal U}^h$, because in this case, the SLA alone can verify the feasibility.

For both tornado lengths, there are decreases in the dislocation by increasing the budget because decision-makers are able to retrofit and recover more locations. For example, increasing the budget from \$0 to \$15M results in around 3000 fewer dislocations in damaged neighborhoods.  Note that this translates in reductions of around 20\% in population dislocation for both tornado lengths by investing \$15M. However, adding another \$15 million will reduce dislocations only by around 500 people because of the great cost of recovery. In fact, the relationship between budget and dislocation exhibits a clear diminishing returns behavior, see Figure~\ref{fig:dimReturns}, where it can be seen that after around \$10M there is a small decrease in population dislocation. From a decision-making point of view, these results implies that it is not reasonable to spend more than \$10M-\$15M in retrofitting and recovery, as the return of investment is not worth it.
\begin{figure}[htb!]
    \centering
    \includegraphics[scale=0.4]{ 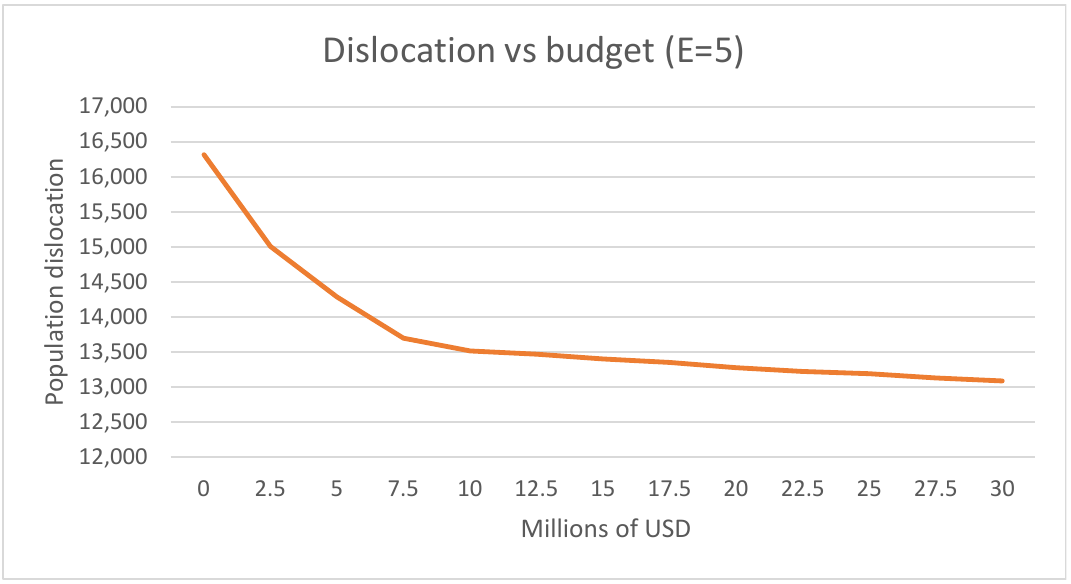}
    \caption{Budget vs population dislocation for a path length of $E=5$ miles.}
    \label{fig:dimReturns}
\end{figure}

More details on the optimal strategies can be found in Figure~\ref{fig:Results-Joplin-Table2}, which maps the locations, their optimal strategies, and the worst-case tornado scenarios for various budgets. It shows that in all experiments the decision-maker uses most of the extra amount of budget from \$15M to \$30M to recover a few locations rather than in retrofit locations that are not in the path of the worst-case tornado. Notice that recovery costs are generally more expensive than retrofitting costs. For example, the recovered location in the {$\infty$-\$15M} map in Figure~\ref{fig:Results-Joplin-Table2} has 470 residents; 382 of which are estimated to be dislocated by spending \$194K for retrofitting. However, the decision-maker has to  spend \$5.5M in recovery to further reduce the estimated dislocation to 295 people for this location. Columns 9-10 in Table~\ref{Tab2:JoplinRobust} show that the proportion of the budget that is assigned to recovery increases by having more budget. In addition, we can observe in Figure~\ref{fig:Results-Joplin-Table2} that the path of worst-case tornadoes does not change by much in different experiments, suggesting that the location of the tornado path is largely driven by the population density in the city center area.

\begin{figure}[!ht]
    \centering
    \captionsetup{justification=centering}
    \includegraphics[width= \textwidth]{ 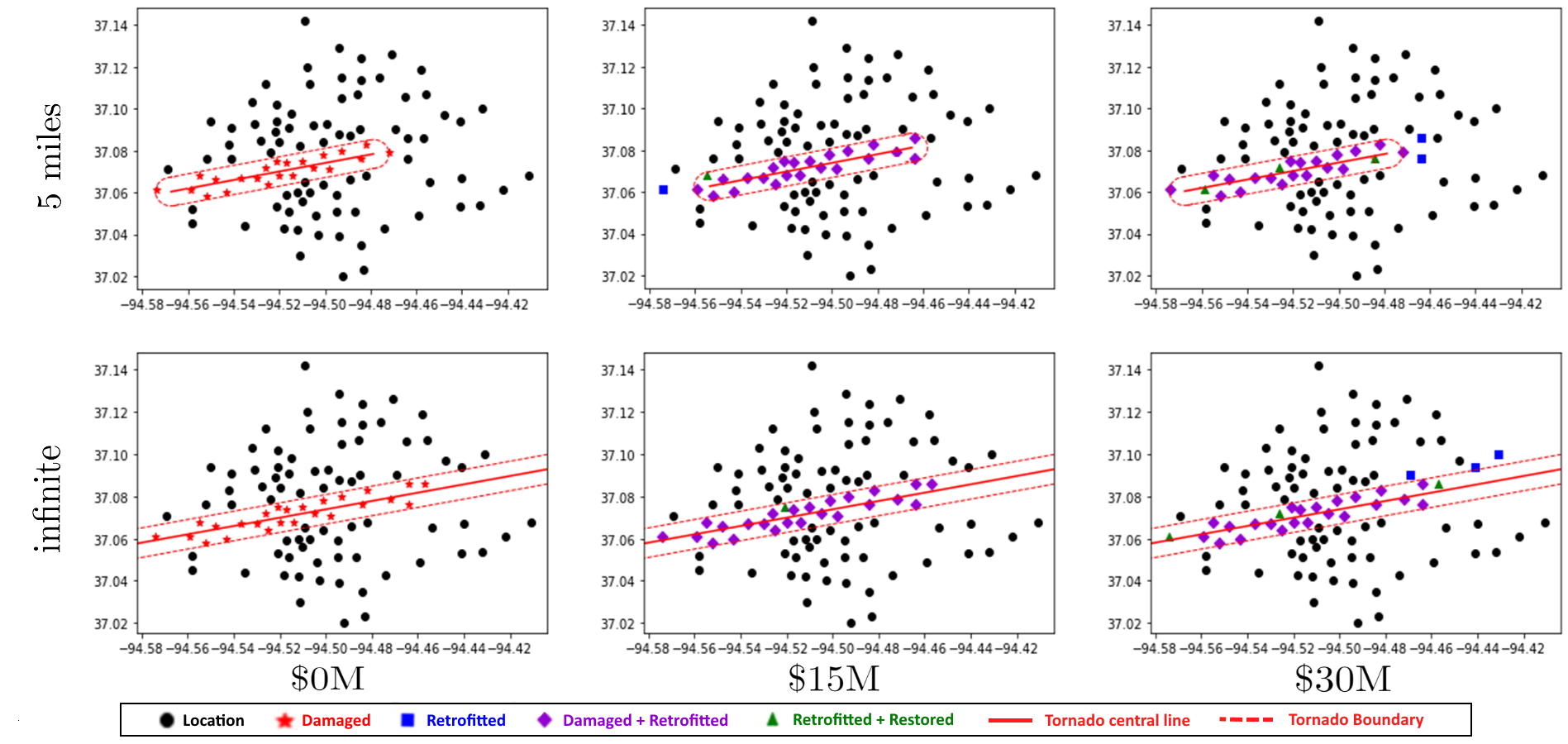}
    \caption{Maps of 100 locations and their retrofitting and recovery strategies under the worst-case tornado scenario for six parameter settings in Table~\ref{Tab2:JoplinRobust}. When retrofitting, the selected strategy is almost exclusively  $R_3$, which is the best, but most expensive, strategy.}
    \label{fig:Results-Joplin-Table2}
\end{figure}

\looseness-1These results suggest, from a worst-case perspective, that it is optimal to invest only in locations in or close to where the worst-case tornado hits and that it is better to invest both in retrofitting and recovery, rather than just in retrofitting. The results are somewhat surprising in the sense that, at least a priori, it would be reasonable to invest all budget in retrofitting across a wider area (particularly because recovery is much more expensive than retrofitting). As such, the results suggest that there is a critical set of clustered locations, in this case the densely populated center of the city, where most of the investments must be concentrated. We suspect, however, that there is a different structure for the optimal actions in data sets with a more homogeneous population distribution.

\subsection{Advantages of robust retrofitting}

In this section we compare the optimal two-stage robust retrofitting and recovery strategies against a policy that makes retrofitting and recovery decisions at random. We assume a $E=5$ miles tornado path and a total budget of $A=$\$15 million. The decision-maker is allowed to spend \$0, \$3, \$9, and \$15 million of the budget in each experiment to randomly retrofit locations with random selections in set $S$; the remaining budget is reserved for the recovery actions. Once a retrofitting strategy $f$ is determined, then the corresponding worst-case second-stage happens (i.e., the value of $\Phi(f)$ is determined, see Equation~\eqref{eq:tornado_bilevel_problem}). Table~\ref{tab:random-retrofit} compares the average, maximum, minimum, and standard deviation of population dislocation for 10 random replications.  Note that for the \$0M case, we do not need multiple replications.

\begin{table}[!ht]
    \centering
    \captionsetup{justification=centering}
     \caption{Comparison of population dislocation statistics for 10 random retrofitting plans with spending a certain amount of budget out of a total \$15M.}
    \begin{adjustbox}{max width=\textwidth}
    \begin{tabular}{c|c|c|c|c}
    \hline
        \multirow{2}{*}{Budget to randomly retrofit} & \multicolumn{4}{c}{Population Dislocation} \\
        & Average & Maximum & Minimum & Standard deviation \\ \hline
        \$0M (0\%) & \multicolumn{3}{c}{16318} & - \\ \hline
        \$3M (20\%) & 16057 & 16263 & 15717 & 163 \\ \hline
        \$9M (60\%) & 15596 & 15905 & 15306 & 196 \\ \hline
        \$15M (100\%) & 15143 & 15451 & 14764 & 189 \\ \hline
        \textbf{Robust optimal retrofitting} & \multicolumn{3}{c}{\textbf{13408}} & - \\ \hline
    \end{tabular}
    \end{adjustbox}
   
    \label{tab:random-retrofit}
\end{table}
\begin{figure}[!ht]
    \centering
    \captionsetup{justification=centering}
    \includegraphics[scale=0.6]{ 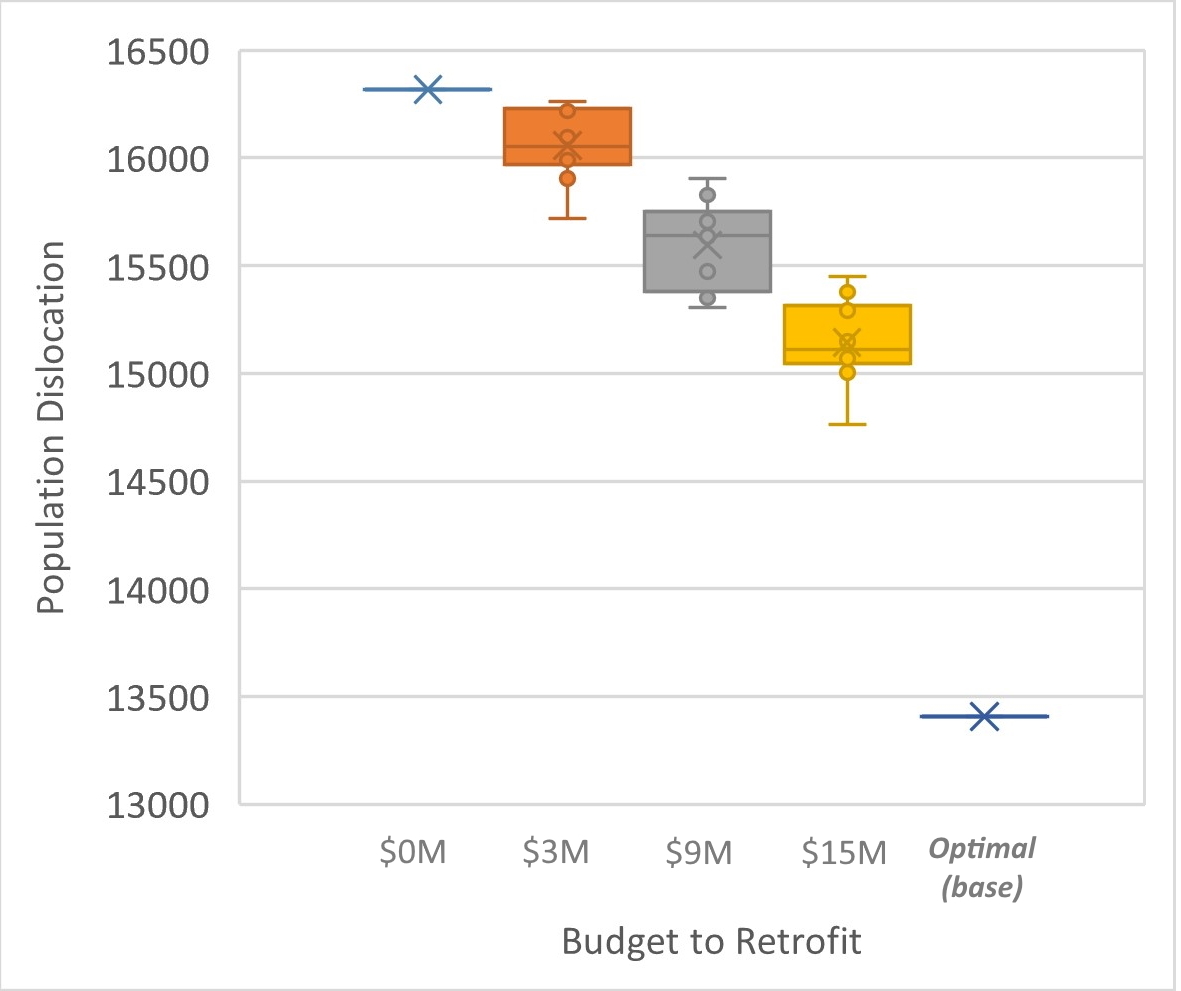}
    \caption{Box and Whisker plots for Table~\ref{tab:random-retrofit} experiments. It shows the reduction trend for population dislocation by considering more budget for retrofitting.}
    \label{fig:random-retrofit}
\end{figure}

The results in Table~\ref{tab:random-retrofit} show that with random policies the average, maximum, and minimum dislocation after the corresponding worst-case tornado decrease  with the percentage of the budget spent in retrofitting. These results are depicted graphically in the Box plots in Figure~\ref{fig:random-retrofit}.
Note that the proposed robust optimization model results in 13408 dislocations, which is noticeable better than the random policies. To further emphasize the value of optimization, recall that the optimal solution with our proposed method in Table~\ref{Tab2:JoplinRobust} spends \$9M out of \$15M (60\% of budget) in retrofitting. By contrast, even the minimum population dislocation with random retrofitting using the same budget proportion is 15306 (a 14\% increase).

We conduct additional simulations to randomly retrofit locations with the full amount of \$15M and \$30M budget as before, but this time the worst-case tornado found in Section~\ref{sec:results Robust strategies in Joplin} is used. The results across 10 replications are shown in Table~\ref{tab:randRetrofit_worstTornado}, which shows that for both cases the robust optimal population dislocation is significantly less than the minimum of population dislocation with random retrofitting. The plot in Figure~\ref{fig:randRetrofit_worstTornado} provides a visualization of the significant gap between robust optimal value (stars) and the distribution of values for random retrofitting plans. We can also observe that population dislocation reduces by increasing the amount of budget. 

\begin{table}[!ht]
    \centering
    \captionsetup{justification=centering}
     \caption{Comparison of population dislocation statistics using ten different  random retrofitting plans and budgets of \$15M and \$30M.}
    \begin{adjustbox}{max width=\textwidth}
    \begin{tabular}{c|c|c|c|c|c}
    \hline
        \multirow{2}{*}{\textbf{Budget}} &\multicolumn{5}{c}{Population dislocation}  \\
         & \textbf{Robust optimal retrofitting} & Average & Maximum & Minimum & Standard deviation \\ \hline
        \textbf{\$15M} & \textbf{13408} & 15314 & 15941 & 14545 & 331 \\ \hline
        \textbf{\$30M} & \textbf{13091} & 14795 & 15248 & 14145 & 294 \\ \hline
    \end{tabular}
    \end{adjustbox}
   
    \label{tab:randRetrofit_worstTornado}
\end{table}
\begin{figure}[!ht]
    \centering
    \captionsetup{justification=centering}
    \includegraphics[scale=0.55]{ 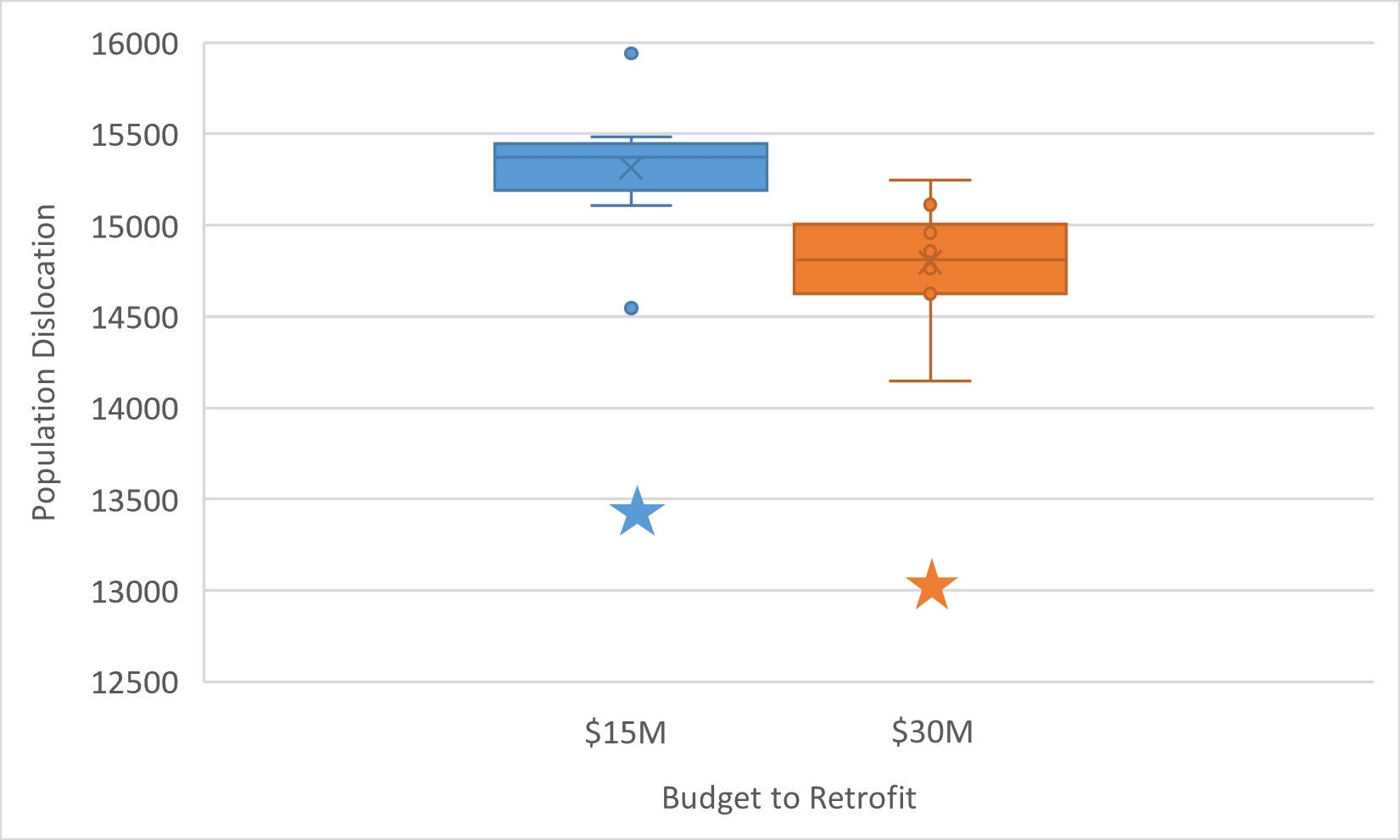}
    \caption{Box and Whisker plots for Table~\ref{tab:randRetrofit_worstTornado} experiments. It shows the gap between the robust optimal plan (stars) and the population dislocation of the random retrofitting strategies for the base worst-case tornado.}
    \label{fig:randRetrofit_worstTornado}
\end{figure}

The results of this section show that deciding using the proposed two-stage robust method yields significant reductions in population dislocation against random policies when evaluating in terms of worst-case outcomes. In the next section we complement these experiments by analyzing how the optimal policy behaves when random tornadoes, instead of the worst-case scenario, happen.

\subsection{Robust strategies in random scenarios}
In this section we fix the optimal retrofitting-recovery plan and compare the robust population dislocation values resulting from having 10 simulated 5-miles tornadoes, see Table~\ref{tab:fixRetrofit_randTronado} and Figure~\ref{fig:fixRetrofit_randTronado}. 

\begin{table}[!ht]
    \centering
    \captionsetup{justification=centering}
    \caption{Comparison of population dislocation statistics for the worst-case tornado and simulated tornadoes by using the optimal retrofitting plan}
    \begin{adjustbox}{max width=\textwidth}
    \begin{tabular}{c|c|c|c|c|c}
    \hline
        \multirow{2}{*}{\textbf{Budget}} & \multicolumn{5}{c}{Population dislocation}  \\
         & \textbf{Worst-case tornado} & Average & Maximum & Minimum & Standard deviation  \\ \hline
        \textbf{\$0M} & \textbf{16318} & 10324 & 13201 & 5233 & 2835 \\ \hline
        \textbf{\$15M} & \textbf{13408} & 9634 & 12720 & 5031 & 2679  \\ \hline
        \textbf{\$30M} & \textbf{13091} & 9416 & 12480 & 4900 & 2607 \\ \hline
    \end{tabular}
    \end{adjustbox}
    
    \label{tab:fixRetrofit_randTronado}
\end{table}
\begin{figure}[!ht]
    \centering
    \captionsetup{justification=centering}
    \includegraphics[scale=0.55]{ 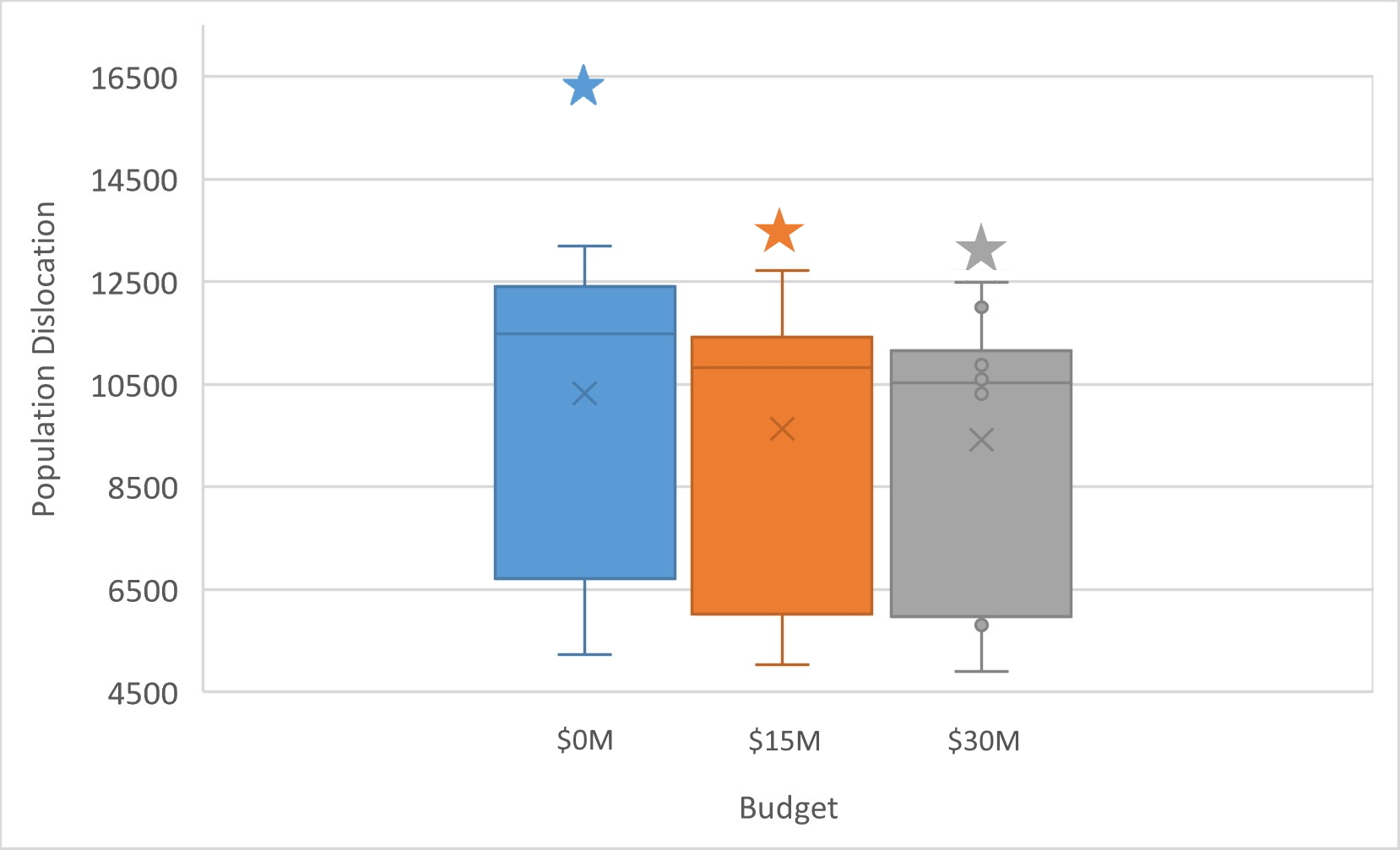}
    \caption{Box and Whisker plots for Table~\ref{tab:fixRetrofit_randTronado} experiments. It shows the distribution of dislocation values for simulated tornadoes and the
    gap to the worst-case tornado.}
    \label{fig:fixRetrofit_randTronado}
\end{figure}

\looseness-1The results for \$0M in Table~\ref{tab:fixRetrofit_randTronado} and Figure~\ref{fig:fixRetrofit_randTronado} show that there is a large gap between the worst-case tornado 16318 (blue star) and the maximum value among simulation results which is 13201. As expected, by having more budget, the average population dislocations in the simulations reduce from 10324 to 9643 for \$15M and to 9416 for \$30M. In addition, the gap between the robust and maximum simulation values for \$15 and \$30M are noticeable smaller than the gap for the no-budget case. 

These observations show the following about the optimal retrofitting-recovery policy: first, that there can be a large variability for population dislocation assuming any tornado is possible. Second, that the variability seems to decrease as more budget is available, and third, that the worst-case scenario is not that far way from simulated scenarios. This last point is further emphasized by the fact that only 10 tornadoes are simulated; one should expect that with more simulated tornadoes the gap between the worst-case and the maximum dislocation observed reduces even more. Further, this last point also shows that the worst-case is not as rare as one might anticipate, and thus that the proposed robust optimization model does not suffer from `over-conservatism,' which is a common drawback of deciding in a robust manner.



\section{Conclusions}\label{sec:conclusion}

Motivated by the need of data-driven methods to allocate resources in retrofitting for tornado-prone areas, we consider a two-stage robust optimization model to determine retrofitting and recovery strategies before and after the realization of an uncertain tornado disaster. We assume that the decision-maker (a combination of public-private agencies) has to allocate a budget in retrofitting and recovery actions, to minimize population dislocation. We propose a two-stage optimization model where the first-stage variables determine the retrofitting strategies and the second-stage variables determine the recovery strategy, assuming that for any retrofitting action, the corresponding worst-case tornado  happens. Given that the uncertainty corresponds to the possible tornado paths, the resulting problem is a two-stage mixed-integer robust optimization problem with a mixed-integer non-linear uncertainty set.

We show that the proposed problem is NP-hard. A standard approach to solve the problem is a column-and-constraint generation algorithm that needs to solve a max-min subproblem at each iteration. We embed a decomposition branch-and-cut algorithm in the column-and-constraint solution method to address the non-convex structure of the bilevel subproblem. We propose two sets of conflict constraints and a custom separation procedure to implement the decomposition branch-and-cut algorithm. Using data from INCORE, we use our model to provide the optimal retrofitting and recovery strategies for the city of Joplin, Missouri against the worst-case scenario.  The results show that the optimal policy exhibits a diminishing returns behavior that sets-in quickly; prioritizes the high-density city center; and prefers to retrofit and recover than to only retrofit. The proposed policy also outperforms other benchmark and does not suffer from `over-conservatism.' 

Future work will focus on modeling tornado damages assuming that there might be several levels of disruptions. Additionally, future research would analyze other social vulnerability indices to measure community resilience under the risk of a tornado. One line of research might study the population dislocation in residential areas with different income levels after certain days of a tornado occurrence, in order to recognize more vulnerable locations that require immediate agencies' attention. Another line of research might take our proposed framework to study other natural disasters like hurricanes with minor modifications in the uncertainty model and the set of retrofitting strategies.

\section*{Acknowledgment}\label{sec:acknowledgment}

This research is partially funded by the National Science Foundation (NSF) (Awards ENG/CMMI \# 2145553 and GEO/RISE \# 2052930) and by the Air Force Office of Scientific Research (AFORS) (Award \# FA9550-22-1-0236). This research was also partially funded by the National Institute of Standards and Technology (NIST) Center of Excellence for Risk-Based Community Resilience Planning through a cooperative agreement with Colorado State University [Awards \# 70NANB20H008 and 70NANB15H044]. The contents expressed in this paper are the views of the authors and do not necessarily represent the opinions or views of the National Institute of Standards and Technology (NIST), the National Science Foundation (NSF), or the Air Force Office of Scientific Research (AFORS).

\bibliography{mybibfile}

\begin{thebibliography}{65}
\providecommand{\natexlab}[1]{#1}
\providecommand{\url}[1]{\texttt{#1}}
\expandafter\ifx\csname urlstyle\endcsname\relax
  \providecommand{\doi}[1]{doi: #1}\else
  \providecommand{\doi}{doi: \begingroup \urlstyle{rm}\Url}\fi

\bibitem[An and Zeng(2015)]{6812211}
Y.~An and B.~Zeng.
\newblock Exploring the modeling capacity of two-stage robust optimization:
  Variants of robust unit commitment model.
\newblock \emph{IEEE Transactions on Power Systems}, 30\penalty0 (1):\penalty0
  109--122, 2015.
\newblock \doi{10.1109/TPWRS.2014.2320880}.

\bibitem[An et~al.(2014)An, Zeng, Zhang, and Zhao]{an2014reliable}
Y.~An, B.~Zeng, Y.~Zhang, and L.~Zhao.
\newblock Reliable p-median facility location problem: two-stage robust models
  and algorithms.
\newblock \emph{Transportation Research Part B: Methodological}, 64:\penalty0
  54--72, 2014.

\bibitem[Ansari(2023)]{github}
M.~Ansari.
\newblock Github repository.
\newblock Online at
  \href{https://github.com/mehdi-ansari/Two_Stage_Robust_Tornado_Problem}{https://github.com/mehdi-ansari/Two\_Stage\_Robust\_Tornado\_Problem},
  Accessed September 2023, September 2023.

\bibitem[Ansari et~al.(2022)Ansari, Borrero, and Lozano]{ansari2022robust}
M.~Ansari, J.~S. Borrero, and L.~Lozano.
\newblock Robust minimum-cost flow problems under multiple ripple effect
  disruptions.
\newblock \emph{INFORMS Journal on Computing}, 2022.

\bibitem[Atamt{\"u}rk and Zhang(2007)]{atamturk2007two}
A.~Atamt{\"u}rk and M.~Zhang.
\newblock Two-stage robust network flow and design under demand uncertainty.
\newblock \emph{Operations Research}, 55\penalty0 (4):\penalty0 662--673, 2007.

\bibitem[Ben-Tal and Nemirovski(1999)]{Ben-TalNemirovski99}
A.~Ben-Tal and A.~Nemirovski.
\newblock Robust solutions to uncertain linear programs.
\newblock \emph{Operations Research Letters}, 25:\penalty0 1--13, 1999.

\bibitem[Ben-Tal et~al.(2004)Ben-Tal, Goryashko, Guslitzer, and
  Nemirovski]{ben2004adjustable}
A.~Ben-Tal, A.~Goryashko, E.~Guslitzer, and A.~Nemirovski.
\newblock Adjustable robust solutions of uncertain linear programs.
\newblock \emph{Mathematical programming}, 99\penalty0 (2):\penalty0 351--376,
  2004.

\bibitem[Ben-Tal et~al.(2011)Ben-Tal, Do~Chung, Mandala, and
  Yao]{ben2011robust}
A.~Ben-Tal, B.~Do~Chung, S.~R. Mandala, and T.~Yao.
\newblock Robust optimization for emergency logistics planning: Risk mitigation
  in humanitarian relief supply chains.
\newblock \emph{Transportation research part B: methodological}, 45\penalty0
  (8):\penalty0 1177--1189, 2011.

\bibitem[Ben-Tal et~al.(2015)Ben-Tal, Hazan, Koren, and Mannor]{ben2015oracle}
A.~Ben-Tal, E.~Hazan, T.~Koren, and S.~Mannor.
\newblock Oracle-based robust optimization via online learning.
\newblock \emph{Operations Research}, 63\penalty0 (3):\penalty0 628--638, 2015.

\bibitem[Bertsimas and Sim(2004)]{BertsimasSim04}
D.~Bertsimas and M.~Sim.
\newblock The price of robustness.
\newblock \emph{Operations Research}, 52\penalty0 (1):\penalty0 35--53, 2004.

\bibitem[Borrero and Lozano(2021)]{borrero2021modeling}
J.~S. Borrero and L.~Lozano.
\newblock Modeling defender-attacker problems as robust linear programs with
  mixed-integer uncertainty sets.
\newblock \emph{INFORMS Journal on Computing}, 33\penalty0 (4):\penalty0
  1570--1589, 2021.

\bibitem[Cheng et~al.(2021)Cheng, Adulyasak, and Rousseau]{cheng2021robust}
C.~Cheng, Y.~Adulyasak, and L.-M. Rousseau.
\newblock Robust facility location under demand uncertainty and facility
  disruptions.
\newblock \emph{Omega}, 103:\penalty0 102429, 2021.

\bibitem[Ding et~al.(2017)Ding, Li, Yang, Jiang, Bie, and Blaabjerg]{7592421}
T.~Ding, C.~Li, Y.~Yang, J.~Jiang, Z.~Bie, and F.~Blaabjerg.
\newblock A two-stage robust optimization for centralized-optimal dispatch of
  photovoltaic inverters in active distribution networks.
\newblock \emph{IEEE Transactions on Sustainable Energy}, 8\penalty0
  (2):\penalty0 744--754, 2017.
\newblock \doi{10.1109/TSTE.2016.2605926}.

\bibitem[Fan and Pang(2019)]{fan2019stochastic}
F.~Fan and W.~Pang.
\newblock Stochastic track model for tornado risk assessment in the us.
\newblock \emph{Frontiers in built environment}, 5:\penalty0 37, 2019.

\bibitem[Farokhnia et~al.(2020)Farokhnia, van~de Lindt, and
  Koliou]{farokhnia2020selection}
K.~Farokhnia, J.~W. van~de Lindt, and M.~Koliou.
\newblock Selection of residential building design requirements to achieve
  community functionality goals under tornado loading.
\newblock \emph{Practice Periodical on Structural Design and Construction},
  25\penalty0 (1):\penalty0 04019035, 2020.

\bibitem[{FEMA}(2003)]{hazus2003multi}
{FEMA}.
\newblock Multi-hazard loss estimation methodology: Earthquake model. {H}azus
  {MH} 2.1.
\newblock 2003.

\bibitem[FEMA(2021)]{FEMA2021}
FEMA.
\newblock Natural hazard retrofit program toolkit.
\newblock
  \href{https://www.fema.gov/sites/default/files/documents/fema_natural-hazards-retrofit-program-tookit.pdf}{https://www.fema.gov/sites/default/files/documents/fema\_natural-hazards-retrofit-program-tookit.pdf},
  Accessed January 2023, February 2021.

\bibitem[FEMA(2022{\natexlab{a}})]{FEMA}
FEMA.
\newblock Tornado--alerts and warnings.
\newblock Online at
  \href{https://community.fema.gov/ProtectiveActions}{https://community.fema.gov/ProtectiveActions},
  Accessed November 2022, November 2022{\natexlab{a}}.

\bibitem[FEMA(2022{\natexlab{b}})]{FEMA2022}
FEMA.
\newblock Resources for repairing, retrofitting and rebuilding after a tornado.
\newblock
  \href{https://www.fema.gov/fact-sheet/resources-repairing-retrofitting-and-rebuilding-after-tornado-0}{https://www.fema.gov/fact-sheet/resources-repairing-retrofitting-and-rebuilding-after-tornado-0},
  Accessed January 2023, March 2022{\natexlab{b}}.

\bibitem[Gabrel et~al.(2014)Gabrel, Lacroix, Murat, and
  Remli]{gabrel2014robust}
V.~Gabrel, M.~Lacroix, C.~Murat, and N.~Remli.
\newblock Robust location transportation problems under uncertain demands.
\newblock \emph{Discrete Applied Mathematics}, 164:\penalty0 100--111, 2014.

\bibitem[Ho-Nguyen and K{\i}l{\i}n{\c{c}}-Karzan(2018)]{ho2018oracles}
N.~Ho-Nguyen and F.~K{\i}l{\i}n{\c{c}}-Karzan.
\newblock Online first-order framework for robust convex optimization.
\newblock \emph{Operations Research}, 66\penalty0 (6):\penalty0 1670--1692,
  2018.

\bibitem[IN-CORE(2022)]{INCORE}
IN-CORE.
\newblock Interdependent networked community resilience modeling environment.
\newblock
  \href{https://incore.ncsa.illinois.edu/}{https://incore.ncsa.illinois.edu/},
  Accessed October 2022, October 2022.

\bibitem[Jabr et~al.(2015)Jabr, Džafić, and Pal]{6824272}
R.~A. Jabr, I.~Džafić, and B.~C. Pal.
\newblock Robust optimization of storage investment on transmission networks.
\newblock \emph{IEEE Transactions on Power Systems}, 30\penalty0 (1):\penalty0
  531--539, 2015.
\newblock \doi{10.1109/TPWRS.2014.2326557}.

\bibitem[Jiang et~al.(2012)Jiang, Zhang, Li, and Guan]{jiang2012benders}
R.~Jiang, M.~Zhang, G.~Li, and Y.~Guan.
\newblock Benders' decomposition for the two-stage security constrained robust
  unit commitment problem.
\newblock In \emph{IIE Annual Conference. Proceedings}, page~1. Institute of
  Industrial and Systems Engineers (IISE), 2012.

\bibitem[Koliou and van~de Lindt(2020)]{koliou2020development}
M.~Koliou and J.~W. van~de Lindt.
\newblock Development of building restoration functions for use in community
  recovery planning to tornadoes.
\newblock \emph{Natural Hazards Review}, 21\penalty0 (2):\penalty0 04020004,
  2020.

\bibitem[Kuligowski et~al.(2014)Kuligowski, Lombardo, Phan, Levitan, and
  Jorgensen]{kuligowski2014final}
E.~D. Kuligowski, F.~T. Lombardo, L.~Phan, M.~L. Levitan, and D.~P. Jorgensen.
\newblock Final report, national institute of standards and technology (nist)
  technical investigation of the may 22, 2011, tornado in joplin, missouri.
\newblock 2014.

\bibitem[Laporte and Louveaux(1993)]{LaporteLouveaux93}
G.~Laporte and F.~V. Louveaux.
\newblock The integer {L}-shaped method for stochastic integer programs with
  complete recourse.
\newblock \emph{Operations Research Letters}, 13:\penalty0 133--142, 1993.

\bibitem[Li et~al.(2021)Li, Zhang, Li, and Wang]{li2021improved}
Y.~Li, F.~Zhang, Y.~Li, and Y.~Wang.
\newblock An improved two-stage robust optimization model for cchp-p2g
  microgrid system considering multi-energy operation under wind power outputs
  uncertainties.
\newblock \emph{Energy}, 223:\penalty0 120048, 2021.

\bibitem[Lin et~al.(2004)Lin, Janak, and Floudas]{lin2004new}
X.~Lin, S.~L. Janak, and C.~A. Floudas.
\newblock A new robust optimization approach for scheduling under uncertainty:
  I. bounded uncertainty.
\newblock \emph{Computers \& chemical engineering}, 28\penalty0 (6-7):\penalty0
  1069--1085, 2004.

\bibitem[Ma et~al.(2018)Ma, Chen, and Wang]{ma2018resilience}
S.~Ma, B.~Chen, and Z.~Wang.
\newblock Resilience enhancement strategy for distribution systems under
  extreme weather events.
\newblock \emph{IEEE Transactions on Smart Grid}, 9\penalty0 (2):\penalty0
  1442--1451, 2018.

\bibitem[MacQueen(1967)]{macqueen1967classification}
J.~MacQueen.
\newblock Classification and analysis of multivariate observations.
\newblock In \emph{5th Berkeley Symp. Math. Statist. Probability}, pages
  281--297, 1967.

\bibitem[Masoomi and van~de Lindt(2018)]{masoomi2018restoration}
H.~Masoomi and J.~W. van~de Lindt.
\newblock Restoration and functionality assessment of a community subjected to
  tornado hazard.
\newblock \emph{Structure and Infrastructure Engineering}, 14\penalty0
  (3):\penalty0 275--291, 2018.

\bibitem[Masoomi et~al.(2018)Masoomi, Ameri, and van~de Lindt]{masoomi2018wind}
H.~Masoomi, M.~R. Ameri, and J.~W. van~de Lindt.
\newblock Wind performance enhancement strategies for residential wood-frame
  buildings.
\newblock \emph{Journal of Performance of Constructed Facilities}, 32\penalty0
  (3):\penalty0 04018024, 2018.

\bibitem[Matthews et~al.(2019)Matthews, Gounaris, and
  Kevrekidis]{matthews2019designing}
L.~R. Matthews, C.~E. Gounaris, and I.~G. Kevrekidis.
\newblock Designing networks with resiliency to edge failures using two-stage
  robust optimization.
\newblock \emph{European Journal of Operational Research}, 279\penalty0
  (3):\penalty0 704--720, 2019.

\bibitem[McAllister et~al.(2015)]{mcallister2015community}
T.~P. McAllister et~al.
\newblock Community resilience planning guide for buildings and infrastructure
  systems, volume i.
\newblock 2015.

\bibitem[Mutapcic and Boyd(2009)]{mutapcic2009oracles}
A.~Mutapcic and S.~Boyd.
\newblock Cutting-set methods for robust convex optimization with pessimizing
  oracles.
\newblock \emph{Optimization Methods \& Software}, 24\penalty0 (3):\penalty0
  381--406, 2009.

\bibitem[Neyshabouri and Berg(2017)]{neyshabouri2017two}
S.~Neyshabouri and B.~P. Berg.
\newblock Two-stage robust optimization approach to elective surgery and
  downstream capacity planning.
\newblock \emph{European Journal of Operational Research}, 260\penalty0
  (1):\penalty0 21--40, 2017.

\bibitem[NOAA(2022{\natexlab{a}})]{NOAA2022}
NOAA.
\newblock National centers for environmental information, tornadoes statistics.
\newblock Online at
  \href{https://www.ncei.noaa.gov/access/monitoring/tornadoes/time-series/12/7?fatalities=true&mean=true}{https://www.ncei.noaa.gov},
  Accessed September 2022, September 2022{\natexlab{a}}.

\bibitem[NOAA(2022{\natexlab{b}})]{SPC22}
NOAA.
\newblock Storm predication center - national weather service.
\newblock Online at
  \href{https://www.spc.noaa.gov/}{https://www.spc.noaa.gov/}, Accessed October
  2022, October 2022{\natexlab{b}}.

\bibitem[NOAA(2022{\natexlab{c}})]{TornadoFAQ}
NOAA.
\newblock Tornadoes frequently asked questions.
\newblock Online at
  \href{https://www.weather.gov/lmk/tornadoesfaq}{https://www.weather.gov/lmk/tornadoesfaq},
  Accessed November 2022, November 2022{\natexlab{c}}.

\bibitem[NOAA(2022{\natexlab{d}})]{supercell}
NOAA.
\newblock Supercells.
\newblock Online at
  \href{https://www.spc.noaa.gov/misc/AbtDerechos/supercells.htm}{https://www.spc.noaa.gov/misc/AbtDerechos/supercells.htm},
  Accessed November 2022, November 2022{\natexlab{d}}.

\bibitem[NOAA(2022{\natexlab{e}})]{tornadocost}
NOAA.
\newblock Weather related fatality and injury statistics.
\newblock Online at
  \href{https://www.weather.gov/hazstat/}{https://www.weather.gov/hazstat/},
  Accessed November 2022, November 2022{\natexlab{e}}.

\bibitem[NOAA(2022{\natexlab{f}})]{weather2022}
NOAA.
\newblock National weather service, the enhanced fujita scale (ef scale).
\newblock \href{https://www.weather.gov/oun/efscale}{https://www.weather.gov},
  Accessed September 2022, September 2022{\natexlab{f}}.

\bibitem[NOAA(2023)]{damAssesment}
NOAA.
\newblock Damage assessment toolkit.
\newblock Online at
  \href{https://apps.dat.noaa.gov/stormdamage/damageviewer/}{https://apps.dat.noaa.gov/stormdamage/damageviewer/},
  Accessed March 2023, March 2023.

\bibitem[NWC(2023)]{okHist}
NWC.
\newblock Historical data of oklahoma county tornadoes.
\newblock Online at
  \href{https://www.weather.gov/oun/tornadodata-county-ok-oklahoma}{https://www.weather.gov/oun/tornadodata-county-ok-oklahoma},
  Accessed January 2023, January 2023.

\bibitem[Ripberger et~al.(2018)Ripberger, Jenkins-Smith, Silva, Czajkowski,
  Kunreuther, and Simmons]{ripberger2018tornado}
J.~T. Ripberger, H.~C. Jenkins-Smith, C.~L. Silva, J.~Czajkowski,
  H.~Kunreuther, and K.~M. Simmons.
\newblock Tornado damage mitigation: Homeowner support for enhanced building
  codes in oklahoma.
\newblock \emph{Risk analysis}, 38\penalty0 (11):\penalty0 2300--2317, 2018.

\bibitem[Saba(2013)]{stabline22}
S.~Saba.
\newblock Line intersecting maximal number of circles (circle ``stabbing''
  problem).
\newblock
  \href{https://sahandsaba.com/line-intersecting-maximal-number-of-circles.html}{https://sahandsaba.com/line-intersecting-maximal-number-of-circles.html},
  Accessed October 2022, October 2013.

\bibitem[Shams et~al.(2021)Shams, Shahabi, MansourLakouraj, Shafie-khah, and
  Catal{\~a}o]{shams2021adjustable}
M.~H. Shams, M.~Shahabi, M.~MansourLakouraj, M.~Shafie-khah, and J.~P.
  Catal{\~a}o.
\newblock Adjustable robust optimization approach for two-stage operation of
  energy hub-based microgrids.
\newblock \emph{Energy}, 222:\penalty0 119894, 2021.

\bibitem[Simmons et~al.(2015)Simmons, Kovacs, and Kopp]{simmons2015tornado}
K.~M. Simmons, P.~Kovacs, and G.~A. Kopp.
\newblock Tornado damage mitigation: Benefit--cost analysis of enhanced
  building codes in oklahoma.
\newblock \emph{Weather, climate, and society}, 7\penalty0 (2):\penalty0
  169--178, 2015.

\bibitem[Standohar-Alfano et~al.(2017)Standohar-Alfano, van~de Lindt, and
  Ellingwood]{standohar2017vertical}
C.~D. Standohar-Alfano, J.~W. van~de Lindt, and B.~R. Ellingwood.
\newblock Vertical load path failure risk analysis of residential wood-frame
  construction in tornadoes.
\newblock \emph{Journal of Structural Engineering}, 143\penalty0 (7):\penalty0
  04017045, 2017.

\bibitem[Stoner and Pang(2021)]{stoner2021tornado}
M.~Stoner and W.~Pang.
\newblock Tornado hazard assessment of residential structures built using
  cross-laminated timber and light-frame wood construction in the us.
\newblock \emph{Natural Hazards Review}, 22\penalty0 (4):\penalty0 04021032,
  2021.

\bibitem[Strader et~al.(2016)Strader, Pingel, and Ashley]{strader2016monte}
S.~M. Strader, T.~J. Pingel, and W.~S. Ashley.
\newblock A monte carlo model for estimating tornado impacts.
\newblock \emph{Meteorological Applications}, 23\penalty0 (2):\penalty0
  269--281, 2016.

\bibitem[Takeda et~al.(2008)Takeda, Taguchi, and
  T{\"u}t{\"u}nc{\"u}]{takeda2008adjustable}
A.~Takeda, S.~Taguchi, and R.~T{\"u}t{\"u}nc{\"u}.
\newblock Adjustable robust optimization models for a nonlinear two-period
  system.
\newblock \emph{Journal of Optimization Theory and Applications}, 136\penalty0
  (2):\penalty0 275--295, 2008.

\bibitem[{The White House}(2022)]{WH2022}
{The White House}.
\newblock Fact sheet: Biden--harris administration launches initiative to
  modernize building codes, improve climate resilience, and reduce energy
  costs.
\newblock
  \href{https://www.whitehouse.gov/briefing-room/statements-releases/2022/06/01/fact-sheet-biden-harris-administration-launches-initiative-to-modernize-building-codes-improve-climate-resilience-and-reduce-energy-costs/}{https://www.whitehouse.gov/briefing-room/statements-releases/2022/06/01/fact-sheet-biden-harris-administration-launches-initiative-to-modernize-building-codes-improve-climate-resilience-and-reduce-energy-costs/},
  June 2022.

\bibitem[Thiele et~al.(2009)Thiele, Terry, and Epelman]{thiele2009robust}
A.~Thiele, T.~Terry, and M.~Epelman.
\newblock Robust linear optimization with recourse.
\newblock \emph{Rapport technique}, pages 4--37, 2009.

\bibitem[van~de Lindt and Dao(2009)]{van2009performance}
J.~W. van~de Lindt and T.~N. Dao.
\newblock Performance-based wind engineering for wood-frame buildings.
\newblock \emph{Journal of Structural Engineering}, 135\penalty0 (2):\penalty0
  169--177, 2009.

\bibitem[Velasquez et~al.(2020)Velasquez, Mayorga, and
  {\"O}zalt{\i}n]{velasquez2020prepositioning}
G.~A. Velasquez, M.~E. Mayorga, and O.~Y. {\"O}zalt{\i}n.
\newblock Prepositioning disaster relief supplies using robust optimization.
\newblock \emph{IISE Transactions}, 52\penalty0 (10):\penalty0 1122--1140,
  2020.

\bibitem[Wang et~al.(2017)Wang, Cao, Pang, and Cao]{wang2017experimental}
J.~Wang, S.~Cao, W.~Pang, and J.~Cao.
\newblock Experimental study on effects of ground roughness on flow
  characteristics of tornado-like vortices.
\newblock \emph{Boundary-Layer Meteorology}, 162:\penalty0 319--339, 2017.

\bibitem[Wang et~al.(2021)Wang, Van De~Lindt, Rosenheim, Cutler, Hartman,
  Sung~Lee, and Calderon]{wang2021effect}
W.~Wang, J.~W. Van De~Lindt, N.~Rosenheim, H.~Cutler, B.~Hartman, J.~Sung~Lee,
  and D.~Calderon.
\newblock Effect of residential building wind retrofits on social and economic
  community-level resilience metrics.
\newblock \emph{Journal of Infrastructure Systems}, 27\penalty0 (4):\penalty0
  04021034, 2021.

\bibitem[Wen(2021)]{wen2021development}
Y.~Wen.
\newblock \emph{Development of Multi-Objective Optimization Model of Community
  Resilience on Mitigation Planning}.
\newblock PhD thesis, University of Oklahoma, 2021.

\bibitem[Xiong et~al.(2017)Xiong, Jirutitijaroen, and
  Singh]{xiong2017distributionally}
P.~Xiong, P.~Jirutitijaroen, and C.~Singh.
\newblock A distributionally robust optimization model for unit commitment
  considering uncertain wind power generation.
\newblock \emph{IEEE Transactions on Power Systems}, 32\penalty0 (1):\penalty0
  39--49, 2017.

\bibitem[Yao et~al.(2009)Yao, Mandala, and Do~Chung]{yao2009evacuation}
T.~Yao, S.~R. Mandala, and B.~Do~Chung.
\newblock Evacuation transportation planning under uncertainty: a robust
  optimization approach.
\newblock \emph{Networks and Spatial Economics}, 9\penalty0 (2):\penalty0 171,
  2009.

\bibitem[Yuan et~al.(2016)Yuan, Wang, Qiu, Chen, Kang, and
  Zeng]{yuan2016robust}
W.~Yuan, J.~Wang, F.~Qiu, C.~Chen, C.~Kang, and B.~Zeng.
\newblock Robust optimization-based resilient distribution network planning
  against natural disasters.
\newblock \emph{IEEE Transactions on Smart Grid}, 7\penalty0 (6):\penalty0
  2817--2826, 2016.

\bibitem[Zeng and Zhao(2013)]{zeng2013solving}
B.~Zeng and L.~Zhao.
\newblock Solving two-stage robust optimization problems using a
  column-and-constraint generation method.
\newblock \emph{Operations Research Letters}, 41\penalty0 (5):\penalty0
  457--461, 2013.

\bibitem[Zhao and Zeng(2012)]{zhao2012robust}
L.~Zhao and B.~Zeng.
\newblock Robust unit commitment problem with demand response and wind energy.
\newblock In \emph{2012 IEEE power and energy society general meeting}, pages
  1--8. IEEE, 2012.

\end{thebibliography}

\newpage

\setcounter{page}{1}
\setcounter{section}{0}

\begin{center}
    \textbf{APPENDIX}
\end{center}

\section{Proof of $NP$-hardness}\label{proof of Proposition Np-hard}
\begin{prop}\label{pr:NP-hardness}
    The optimization problem~\eqref{eq:two_stage_model} is NP-hard.
\end{prop}
\noindent\emph{Proof:} Consider an arbitrary instance of the knapsack problem with $N$ items, weights $b_i, \forall i \in N$, values $q_i, \forall i \in N$, capacity $B$, and value $Q$. 
We build an instance of problem \eqref{eq:two_stage_model} that is equivalent to the knapsack problem. 
Consider an instance of problem~\eqref{eq:two_stage_model} with $|S|=1$ and $|P|=2$. Observe that because $|S|=1$, $f$ has to be equal to a vector of ones. 
Suppose that the coordinates of the locations and the values of $E$ and $\Delta$ are such that $z_\ell=1$ for all $\ell\in L$ is feasible in $\mathcal{U}$. For simplicity, we drop index $s$ from the notation. Also, note that $r_{\ell 2}$ can be replaced by $1 - r_{\ell 1}$ for all $\ell \in L$ because only two recovery strategies are assumed. 
The specific case of the problem can be rewritten as
    \begin{subequations}
\label{eq:second_stage_proof}
\begin{align}
    v = \min & \sum_{\ell\in L} (g_{\ell 1} - g_{\ell 2}) r_{\ell 1} + \sum_{\ell\in L} (w_{\ell} + g_{\ell 2})  \label{pr_eq:ssp0}
    \\
    \text{s.t. }& \sum_{\ell \in L}(c_{\ell 1} - c_{\ell 2})r_{\ell 1} \le A - \sum_{\ell \in L} (d_\ell + c_{\ell 2})\label{pr_eq:ssp1}
    \\
    &r_{\ell 1}\in\{0,1\}^{|L|}.
\end{align}
\end{subequations}

Consider the example of problem~\eqref{eq:second_stage_proof} where $L=N$, $c_{\ell 1}- c_{\ell 2} = b_{\ell}$ and $g_{\ell 2}- g_{\ell 1} = q_{\ell}$ for each item $\ell \in L$, and $A - \sum_{\ell \in L} (d_\ell + c_{\ell 2})=B$. The answer of knapsack problem for some instance $L^* \subseteq N$ is YES (i.e., $\sum_{i \in L^*} q_i \geq Q$), if and only if ${v}\le \sum_{\ell\in L} (w_{\ell} + g_{\ell 2}) - Q$.\epf

\section{Proof of Proposition \ref{pr:inf_pair}}
\label{proof of Proposition1}
\noindent\emph{Proof:} Suppose $(\ell_1,\ell_2) \in \Omega$, and for the sake of contradiction that $z_{\ell_1}=z_{\ell_2}=1$. For the constraint in \eqref{eq:uncertain_LS_coverage}, we have
\begin{subequations}
\begin{align}
    \|(1-t_{\ell_1})e_0+t_{\ell_1}e_{1}-(x_{\ell_1},y_{\ell_1})\| \le \Delta
    \\
    \|(1-t_{\ell_2})e_0+t_{\ell_2}e_{1}-(x_{\ell_2},y_{\ell_2})\| \le \Delta.
\end{align}
\end{subequations}
The sum of two above inequalities results
\begin{equation}
    \|(1-t_{\ell_1})e_0+t_{\ell_1}e_{1}-(x_{\ell_1},y_{\ell_1})\| + \|(x_{\ell_2},y_{\ell_2}) - (1-t_{\ell_2})e_0 - t_{\ell_2}e_{1}\| \le 2\Delta,
\end{equation}
and by applying the triangle inequality $\|X+Y\|\le\|X\|+\|Y\|$, we conclude that
\begin{equation}
    \|(x_{\ell_2}-x_{\ell_1}, y_{\ell_2}- y_{\ell_1}) + (t_{\ell_2}-t_{\ell_1})(e_0-e_1)\| \le 2\Delta.
\end{equation}
Finally, we employ the fact that $\|X\|-\|Y\|\le\|X-Y\|$ to reach
\begin{equation}
\label{compact_expans}
    \|(x_{\ell_2}-x_{\ell_1}, y_{\ell_2}- y_{\ell_1})\| - \|(t_{\ell_1}-t_{\ell_2})(e_0-e_1)\| \le 2\Delta.
\end{equation}
Now, we show \eqref{compact_expans} cannot be true because $(\ell_1,\ell_2)\in\Omega$, i.e.,
\begin{equation}
\label{eq:proof1_condOmega}
     \|(x_{\ell_2}-x_{\ell_1}, y_{\ell_2}- y_{\ell_1})\|>2\Delta+E.
\end{equation}
By constraint \eqref{eq:uncertain_LS_length} and using $\|aY\|=|a|\|Y\|$, we have
\begin{equation} \label{first_term}
    \|(t_{\ell_1}-t_{\ell_2})(e_0-e_1)\| = |t_{\ell_1}-t_{\ell_2}|\|e_0-e_1\| \le |t_1-t_2|E \le E.
\end{equation}
Note that $|t_1-t_2|\le1$. So, the two above inequalities imply that
\begin{equation}
\begin{aligned}
\label{contradiction_res}
    &\begin{cases}
    \|(x_{\ell_2}-x_{\ell_1}, y_{\ell_2}- y_{\ell_1})\|>2\Delta+E
    \\
    -\|(t_1-t_2)(e_0-e_1)\|_1\ge -E
    \end{cases}
    \\
    & \quad \implies \|(x_{\ell_2}-x_{\ell_1}, y_{\ell_2}- y_{\ell_1})\| - \|(t_{\ell_1}-t_{\ell_2})(e_0-e_1)\| > 2\Delta.
\end{aligned}
\end{equation}
Observe that \eqref{compact_expans} violates the true inequality \eqref{contradiction_res}. Therefore, the cut $z_{\ell_1}+z_{\ell_2}\le1$ is valid because $z_{\ell_1}=z_{\ell_2}=1$ is an infeasible solution under the assumption $(\ell_1,\ell_2)\in\Omega$.\epf

\section{Proof of Proposition \ref{pr:triple_validity}}
\label{proof of Proposition2}

This section provides a more detailed discussion (and a more detailed proof) for Proposition~\ref{pr:triple_validity}. 
\noindent\emph{Proof:} We are interested to identify the division of space $R$ into feasible ($R_{\ell_1,\ell_2}'$) and infeasible ($R_{\ell_1,\ell_2}''$) subspaces for two locations $\ell_1$ and $\ell_2$ with coordinates $(x_{\ell_1},y_{\ell_1})$ and $(x_{\ell_2},y_{\ell_2})$, respectively. For the sake of simplicity, we assume that $x_{\ell_1} \le x_{\ell_2}$ throughout the proof; the results  can be easily extended to the case $x_{\ell_1} > x_{\ell_2}$.

\begin{lemma}\label{lemma:boundaries}
    Lines $\lambda_1,\ldots,\lambda_6$ in \eqref{eq:P1}-\eqref{eq:P6} determine boundaries of subspaces.
\end{lemma}

\noindent\emph{Proof:} 
We find feasible slope $m$ and y-intercept $q$ of lines in the form of $y=mx+q$ with a distance $\Delta$ from both locations $\ell_1$ and $\ell_2$ which means the following equalities must simultaneously hold,
\begin{subequations}
\label{eq:validproof_dist-delta}
\begin{align}
    \frac{|y_{\ell_1}-mx_{\ell_1}-q|}{\sqrt{1+m^2}} = \Delta, \label{eq:validproof_dist-delta1}
    \\
    \frac{|y_{\ell_2}-mx_{\ell_2}-q|}{\sqrt{1+m^2}} = \Delta. \label{eq:validproof_dist-delta2}
\end{align}
\end{subequations}
The above equalities imply two possible cases:
\begin{equation}
\begin{aligned}
    \frac{|y_{\ell_1}-mx_{\ell_1}-q|}{\sqrt{1+m^2}} &= \frac{|y_{\ell_2}-mx_{\ell_2}-q|}{\sqrt{1+m^2}} \\
    \implies
    &\begin{cases}
    y_{\ell_1}-mx_{\ell_1}= y_{\ell_2}-mx_{\ell_2} \ (\textbf{Case 1})
    \\
    y_{\ell_1}-mx_{\ell_1}-q = mx_{\ell_2}+q - y_{\ell_2} \ (\textbf{Case 2})
    \end{cases}.
\label{eq:validproof-equal}
\end{aligned}
\end{equation}
\textbf{Case 1.}
\begin{equation}
    m = \frac{y_{\ell_2}- y_{\ell_1}}{x_{\ell_2}-x_{\ell_1}} = \tan(\theta) \ (\text{by the definition}).
\end{equation}
By replacing $m=\tan(\theta)$ in \eqref{eq:validproof_dist-delta1}, we will have
\begin{subequations}
    \begin{align}
        &\frac{|y_{\ell_1}-\tan(\theta) x_{\ell_1}-q|}{\sqrt{1+\tan^2(\theta)}} = \Delta &
        \\
        \implies & |y_{\ell_1}-\tan(\theta) x_{\ell_1}-q||\cos(\theta)| = \Delta & (\sqrt{1+\tan^2(\theta)}= \frac{1}{|\cos(\theta)|})
        \\
        \implies & |y_{\ell_1}-\tan(\theta) x_{\ell_1}-q| = \frac{\Delta}{\cos(\theta)} & (x_{\ell_1} \le x_{\ell_2} \implies \cos(\theta) \ge 0)
        \\
        \implies &
        \begin{cases}
        q' = y_{\ell_1} - \tan(\theta) x_{\ell_1} - \frac{\Delta}{\cos(\theta)}
        \\
        q''= y_{\ell_1} - \tan(\theta) x_{\ell_1} + \frac{\Delta}{\cos(\theta)}
        \end{cases}
    \end{align}
\end{subequations}
So, the tornado paths with formulations 
\begin{subequations}
\begin{align}
    y= \tan(\theta) (x-x_{\ell_1}) + y_{\ell_1} - \frac{\Delta}{\cos(\theta)},
    \label{eq:validproof-line1}
    \\
    y= \tan(\theta) (x-x_{\ell_1}) + y_{\ell_1} + \frac{\Delta}{\cos(\theta)},
    \label{eq:validproof-line2}
\end{align}
\end{subequations}
are two boundaries that cover both locations (observe that any other $q$ between $q'$ and $q''$ ($q'\le q\le q''$) corresponding to the slope $m=\tan\theta$ forms a tornado path that has \textbf{a distance at most $\Delta$} to the locations $\ell_1$ and $\ell_2$). The tornado in \eqref{eq:validproof-line1} covers all locations in distance $\Delta$ between the lines $y=\tan\theta (x- x_{\ell_1}) + y_{\ell_1} - \frac{2\Delta}{\cos\theta}$ (see \eqref{eq:P2} for line $\lambda_2$) and $y=\tan\theta (x- x_{\ell_1}) + y_{\ell_1}$. Similarly, the tornado \eqref{eq:validproof-line2} covers all points between $y=\tan\theta (x- x_{\ell_1}) + y_{\ell_1}$ and $y=\tan\theta (x- x_{\ell_1}) + y_{\ell_1} + \frac{2\Delta}{\cos\theta}$ (see \eqref{eq:P1} for line $\lambda_1$).

Replacing $m=\tan\theta$ in \eqref{eq:validproof_dist-delta2} will give the same boundaries as \eqref{eq:validproof_dist-delta1}.
\\
\textbf{Case 2.}
\begin{equation}
    q = \frac{y_{\ell_1}+y_{\ell_2} - m(x_{\ell_1}+x_{\ell_2})}{2}.
\end{equation}
By replacing $q = \frac{y_{\ell_1}+y_{\ell_2} - m(x_{\ell_1}+x_{\ell_2})}{2}$ in \eqref{eq:validproof_dist-delta1}, we have
\begin{equation}
\label{eq:validproof_case2}
    \begin{aligned}
        &\frac{|y_{\ell_1}-mx_{\ell_1}-\frac{y_{\ell_1}+y_{\ell_2} - m(x_{\ell_1}+x_{\ell_2})}{2}|}{\sqrt{1+m^2}} = \Delta
        \\
        \implies & 
        \frac{|y_{\ell_2}- y_{\ell_1}-m(x_{\ell_2}-x_{\ell_1})|}{\sqrt{1+m^2}} = 2\Delta
    \end{aligned}
\end{equation}
Instead of solving \eqref{eq:validproof_case2} for $m$ which leads to a complicated quadratic equation, we exploit the form of equality to compute possible values of $m$. Observe that \eqref{eq:validproof_case2} can be interpreted as a line passing through the origin and has a distance $2\Delta$ from the point $(x_{\ell_2}-x_{\ell_1} , y_{\ell_2}- y_{\ell_1})$  (see Figure \ref{fig:validineq_proof1.jpg}). 
\begin{figure}[h!]
    \centering
    \includegraphics[width=8cm]{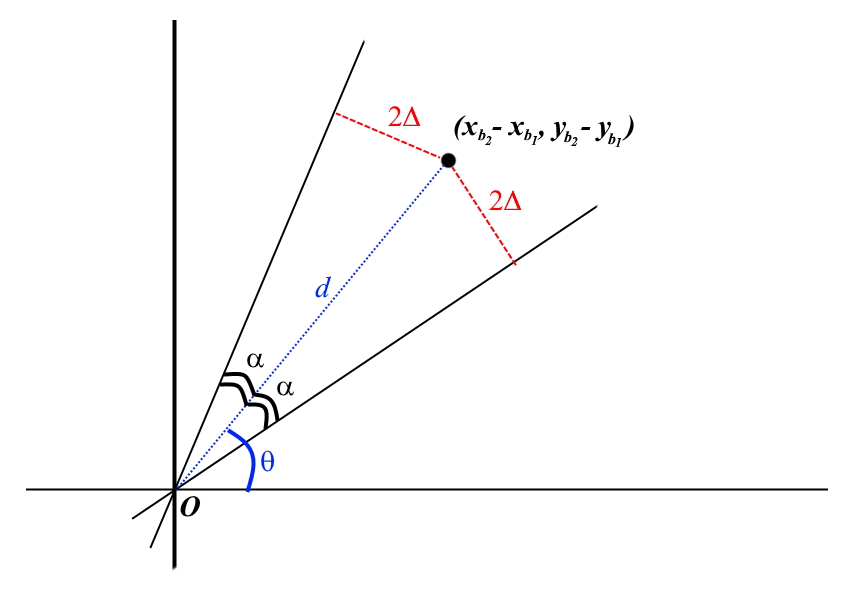}
    \caption{Two lines passing the origin and have a distance $2\Delta$ from the point $(x_{\ell_2}-x_{\ell_1} , y_{\ell_2}- y_{\ell_1})$}
    \label{fig:validineq_proof1.jpg}
\end{figure}
Now, we discuss that the angles corresponding to the possible slopes of $m$ must be 
\begin{align}
    m' = \tan{\left(\arctan{\frac{ y_{\ell_2}- y_{\ell_1}}{ x_{\ell_2}- x_{\ell_1}}} + \arcsin{\frac{2\Delta}{d}}\right)},
    \\
    m'' = \tan{\left(\arctan{\frac{ y_{\ell_2}- y_{\ell_1}}{ x_{\ell_2}- x_{\ell_1}}} - \arcsin{\frac{2\Delta}{d}}\right)},
\end{align}
where $\arctan{\frac{ y_{\ell_2}- y_{\ell_1}}{ x_{\ell_2}- x_{\ell_1}}} = \theta$ and $\arcsin{\frac{2\Delta}{D(\ell_1,\ell_2)}} = \alpha$ by the definition. Therefore, 
\begin{align}
    &m' = \tan{(\theta+\alpha)},
    \\
    &m'' = \tan{(\theta-\alpha)}.
\end{align}
So, the resulting tornado paths are 
\begin{equation}
\begin{aligned}
    y &= \tan(\theta+\alpha)x + \frac{y_{\ell_1}+y_{\ell_2} - \tan(\theta+\alpha)(x_{\ell_1}+x_{\ell_2})}{2} 
    \\
    &=  \tan(\theta+\alpha) (x - \frac{x_{\ell_1}+x_{\ell_2}}{2}) +  \frac{y_{\ell_1}+y_{\ell_2}}{2} \label{eq:validproof_rotatedline1}
\end{aligned} 
\end{equation}
\begin{equation}
\begin{aligned}
    y &= \tan(\theta-\alpha)x + \frac{y_{\ell_1}+y_{\ell_2} - \tan(\theta-\alpha)(x_{\ell_1}+x_{\ell_2})}{2} 
    \\ &=  \tan(\theta-\alpha) (x - \frac{x_{\ell_1}+x_{\ell_2}}{2}) +  \frac{y_{\ell_1}+y_{\ell_2}}{2} \label{eq:validproof_rotatedline2}
\end{aligned} 
\end{equation}
Observe that these tornado paths pass through the midpoint $(\frac{x_{\ell_1}+x_{\ell_2}}{2}, \frac{y_{\ell_1}+y_{\ell_2}}{2})$. So, two edges of the tornado path in \eqref{eq:validproof_rotatedline1} with distance $\Delta$ will be $y=\tan(\theta+\alpha) (x - x_{\ell_1}) +  y_{\ell_1}$ and $y=\tan(\theta+\alpha) (x - x_{\ell_2}) +  y_{\ell_2}$ (see lines $\lambda_3$ and $\lambda_4$ in \eqref{eq:P3} and \eqref{eq:P4}). Likewise, the $\Delta$ edges of \eqref{eq:validproof_rotatedline2} are $y=\tan(\theta-\alpha) (x - x_{\ell_1}) +  y_{\ell_1}$ and $y=\tan(\theta-\alpha) (x - x_{\ell_2}) +  y_{\ell_2}$ (see lines $\lambda_5$ and $\lambda_6$ in \eqref{eq:P5} and \eqref{eq:P6}). Replacing $q$ in \eqref{eq:validproof_dist-delta2} leads to the same equality in \eqref{eq:validproof_case2}, and so the same conclusions.\epf

In the following, we determine the feasible pair $(\beta,q(\beta))$ for the line $y=\tan(\beta) x+q(\beta)$ to cover both locations $\ell_1$ and $\ell_2$ within the distance $\Delta$, where $\beta=\arctan(m)$ is the angle of the slope of the line and $q(\beta)$ is the $y$-intercept of the line given angle $\beta$. It can be shown that any feasible slope angle $\beta$ must be inside the interval $[\theta-\alpha,\theta+\alpha]$ as follows.
\begin{subequations}
\label{eq:validproof_dist-le-delta1}
\begin{align}
    \frac{|y_{\ell_1}-\tan(\beta) x_{\ell_1}-q(\beta)|}{\sqrt{1+\tan^2\beta}} \le \Delta, \label{eq:validproof_dist-le-delta1x}
    \\
    \frac{|y_{\ell_2}-\tan(\beta) x_{\ell_2}-q(\beta)|}{\sqrt{1+\tan^2\beta}} \le \Delta. \label{eq:validproof_dist-le-delta2x}
\end{align}
\end{subequations}
The summation of \eqref{eq:validproof_dist-le-delta1x} and \eqref{eq:validproof_dist-le-delta2x} results in
\begin{equation}
     \frac{|-y_{\ell_1}+\tan(\beta) x_{\ell_1}+q(\beta)|+|y_{\ell_2}-\tan(\beta) x_{\ell_2}-q(\beta)|}{\sqrt{1+\tan^2\beta}} \le 2\Delta
\end{equation}
By the triangle inequality $\|X+Y\|\le\|X\|+\|Y\|$, we conclude
\begin{equation}\label{pr:distanceSUM}
    \frac{|y_{\ell_2}-y_{\ell_1}-\tan(\beta) (x_{\ell_2}-x_{\ell_1})|}{\sqrt{1+\tan^2\beta}} \le 2\Delta.
\end{equation}
Based on Case 2 in the proof of Lemma~\ref{lemma:boundaries}, \eqref{pr:distanceSUM} implies that $\beta\in [\theta-\alpha,\theta+\alpha]$. 

Now, we compute the feasible $q(\beta)$ for given $\beta \in [\theta-\alpha,\theta+\alpha]$. From \eqref{eq:validproof_dist-le-delta1x} and \eqref{eq:validproof_dist-le-delta2x}, we have 
\begin{subequations}
\begin{align}
    &y_{\ell_1}-\tan(\beta) x_{\ell_1} -\frac{\Delta}{|\cos(\beta)|} \le q(\beta) \le y_{\ell_1}-\tan(\beta) x_{\ell_1} +\frac{\Delta}{|\cos(\beta)|} \label{eq:pr-yIntercept-ineq1x}
    \\
    &y_{\ell_2}-\tan(\beta) x_{\ell_2} -\frac{\Delta}{|\cos(\beta)|} \le q(\beta) \le y_{\ell_2}-\tan(\beta) x_{\ell_2} +\frac{\Delta}{|\cos(\beta)|}. \label{eq:pr-yIntercept-ineq2x}
\end{align}
\end{subequations}
Note that $\sqrt{1+\tan^2\beta}=1/{|\cos(\beta)|}$. For the sake of simplicity, we assume hereafter that $\theta+\alpha< \pi/2$ and $\theta-\alpha > -\pi/2$ (which implies $\cos(\beta) > 0$) and then discuss two cases for a feasible $q(\beta)$ followed by the inequalities \eqref{eq:pr-yIntercept-ineq1x} and \eqref{eq:pr-yIntercept-ineq2x}:
\begin{enumerate}
    \item $\theta \le \beta \le \theta + \alpha \implies \tan(\theta) \le \tan(\beta) \implies \frac{y_{\ell_2}-y_{\ell_1}}{x_{\ell_2}-x_{\ell_1}} \le \tan(\beta)$:
    \begin{equation}
    \label{pr:y-inter interval1}
        y_{\ell_1}-\tan(\beta) x_{\ell_1} -\frac{\Delta}{\cos(\beta)} \le q(\beta) \le y_{\ell_2}-\tan(\beta) x_{\ell_2} +\frac{\Delta}{\cos(\beta)}.
    \end{equation}

    \item $\theta - \alpha \le \beta \le \theta \implies \tan(\theta) \ge \tan(\beta) \implies \frac{y_{\ell_2}-y_{\ell_1}}{x_{\ell_2}-x_{\ell_1}} \ge \tan(\beta)$:
    \begin{equation}
    \label{pr:y-inter interval2}
        y_{\ell_2}-\tan(\beta) x_{\ell_2} -\frac{\Delta}{\cos(\beta)} \le q(\beta) \le y_{\ell_1}-\tan(\beta) x_{\ell_1} +\frac{\Delta}{\cos(\beta)}.
    \end{equation}
\end{enumerate}
Figure \ref{fig:possible y-intercept.png} visualizes an example of the feasible ranges~\eqref{pr:y-inter interval1} and~\eqref{pr:y-inter interval2} for $q(\beta)$ given angle $\beta \in [\theta-\alpha, \theta+\alpha]$. 

\begin{figure}[h!]
    \centering
    \includegraphics[width = 6cm]{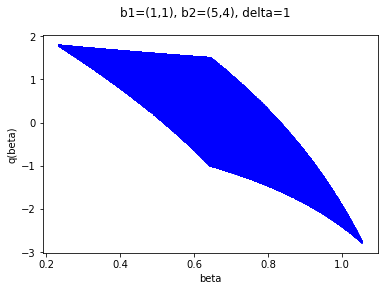}
    \caption{The admissible $y$-intercept $q(\beta)$ according to slope angle $\beta$ for two locations located at $(1,1)$ and $(5,4)$ with given $\Delta=1$. Note that $\theta=\arctan(3/4)$ and $\alpha=\arcsin(2/5)$.}
    \label{fig:possible y-intercept.png}
\end{figure}

Inequalities \eqref{pr:y-inter interval1} and \eqref{pr:y-inter interval2} imply line $y=\tan(\beta) x+q(\beta)$ for given $\beta \in [\theta-\alpha, \theta+\alpha]$ is \emph{feasible} (covers $\ell_1$ and $\ell_2$ within a distance of $\Delta$), if and only if a point $(\bar x, \bar y)$ on the line  meets
\begin{equation}
\label{pr:feasible_lines}
\begin{cases}
     \beta\in [\theta,\theta+\alpha]:\bar y \in [\tan(\beta) (\bar x -x_{\ell_1}) + y_{\ell_1} -\frac{\Delta}{\cos(\beta)}, \tan(\beta) (\bar x- x_{\ell_2}) + y_{\ell_2} +\frac{\Delta}{\cos(\beta)}],
    \\
    \beta\in [\theta-\alpha,\theta]:\bar y \in [\tan(\beta) (\bar x-x_{\ell_2}) + y_{\ell_2} -\frac{\Delta}{\cos(\beta)}, \tan(\beta) (\bar x- x_{\ell_1}) + y_{\ell_1} +\frac{\Delta}{\cos(\beta)}].
    \end{cases}
\end{equation}

\begin{remark}
Given the angle $\beta \in [\theta-\alpha, \theta+\alpha]$, a feasible line $y=\tan(\beta) x+q(\beta)$ covers a point $(x_0,y_0)$ within the distance $\Delta$, {if and only if} 
\begin{equation}
\begin{cases}
    \beta\in [\theta,\theta+\alpha]: y_0 \in [\tan(\beta) (x_0-x_{\ell_1}) + y_{\ell_1} -\frac{2\Delta}{\cos(\beta)}, \tan(\beta) (x_0- x_{\ell_2}) + y_{\ell_2} +\frac{2\Delta}{\cos(\beta)}],
    \\
    \beta\in [\theta-\alpha,\theta]:  y_0 \in [\tan(\beta) (x_0-x_{\ell_2}) + y_{\ell_2} -\frac{2\Delta}{\cos(\beta)}, \tan(\beta) (x_0- x_{\ell_1}) + y_{\ell_1} +\frac{2\Delta}{\cos(\beta)}].
    \end{cases}
\end{equation}
\label{remark:coverage} 
\end{remark}

Regarding Remark~\ref{remark:coverage}, to obtain an inclusive set of feasible coverage for any line $y=\tan(\beta) x+q(\beta)$ where $\beta\in [\theta-\alpha, \theta+\alpha]$, define $$I(x) = \bigcup_{\theta\le\beta\le\theta+\alpha} [\tan(\beta) (x-x_{\ell_1}) + y_{\ell_1} -\frac{2\Delta}{\cos(\beta)}, \tan(\beta) (x - x_{\ell_2}) + y_{\ell_2} +\frac{2\Delta}{\cos(\beta)}],$$ and $$J(x) = \bigcup_{\theta-\alpha\le\beta\le\theta}[\tan(\beta) (x-x_{\ell_2}) + y_{\ell_2} -\frac{2\Delta}{\cos(\beta)}, \tan(\beta) (x - x_{\ell_1}) + y_{\ell_1} +\frac{2\Delta}{\cos(\beta)}].$$
A point $(x,y)$ is covered within distance $\Delta$ of \textit{some} feasible line if and only if $y \in I(x)\cup J(x)$. 

In the following, we find the upper and lower bounds of the set $I(x)\cup J(x)$ to create the set of infeasible $\Gamma_{\ell_1, \ell_2}$ for some location $b_3 \in L$ which $y_{b_3} \notin I(x_{b_3})\cup J(x_{b_3}), \forall \beta\in [\theta-\alpha,\theta+\alpha]$. First, the upper bound of set $I(x)$ is discussed as follows. 

Define
$$
x^I_u = x_{\ell_2} - 2\Delta \left( \frac{\sqrt{1+\tan^2(\theta+\alpha)} - \sqrt{1+\tan^2(\theta)}}{\tan(\theta+\alpha) - \tan(\theta)} \right),
$$
which is the intersection of two extreme lines $\lambda_1$ and $\lambda_3$ (see Lemma~\ref{lemma:boundaries} and note that $y=\tan\theta (x- x_{\ell_2}) + y_{\ell_2} + \frac{2\Delta}{\cos\theta}$ and $y=\tan(\theta+\alpha) (x - x_{\ell_2}) +  y_{\ell_2} + \frac{2\Delta}{\cos (\theta+\alpha)} $ are other representations of $\lambda_1$ and $\lambda_3$). 
We claim that line $\lambda_1$ for $x \le x^I_u$ and line $\lambda_3$ for $x \ge x^I_u$ are the upper bounds of $I(x)$.

\begin{lemma}
\label{lemma:increasing Function}
Function $f(z)=\frac{\sqrt{1+\tan^2 z} - \sqrt{1+\tan^2 a}}{\tan z - \tan a}$ is strictly increasing.
\end{lemma}
\noindent\emph{Proof:} 
We show the derivative of $f(z)$ is strictly positive:
\begin{equation}
\begin{aligned}
   f'(z) =& \frac{\frac{\tan z \sec^2z}{\sqrt{1+\tan^2 z}}(\tan z - \tan a)- \sec^2z(\sqrt{1+\tan^2 z} - \sqrt{1+\tan^2 a})}{(\tan z - \tan a)^2}
\\
= & \frac{\sec^2z}{(\sqrt{1+\tan^2 z})(\tan z - \tan a)^2} \\ & \ \times (\tan^2z - \tan z \tan a - 1 -\tan^2z + \sqrt{1+ \tan^2z + \tan^2a + \tan^2z \tan^2a})
\\
= & \frac{\sec^2z}{(\sqrt{1+\tan^2 z})(\tan z - \tan a)^2} \\ & \ \times (\sqrt{(1+\tan z \tan a)^2 + \tan^2z + \tan^2a - 2 \tan z \tan a} - \tan z \tan a - 1 )
\\
= &  \frac{\sec^2z}{(\sqrt{1+\tan^2 z})(\tan z - \tan a)^2} \\ & \ \times (\sqrt{(1+\tan z \tan a)^2 + (\tan z - \tan a)^2} - (\tan z \tan a + 1)) > 0.
\end{aligned}
\end{equation}
\epf

Using Lemma~\ref{lemma:increasing Function}, we conclude:
\begin{enumerate}
    \item $x \le x^I_u$: 
    \begin{equation}
        \begin{aligned}
            x&\le x_{\ell_2} - 2\Delta \left( \frac{\sqrt{1+\tan^2(\theta+\alpha)} - \sqrt{1+\tan^2(\theta)}}{\tan(\theta+\alpha) - \tan(\theta)} \right)
            \\
            & \le x_{\ell_2} - 2\Delta \left( \frac{\sqrt{1+\tan^2(\beta)} - \sqrt{1+\tan^2(\theta)}}{\tan(\beta) - \tan(\theta)} \right) & \forall \beta \in [\theta, \theta + \alpha]
             \\ \implies
            & \tan(\theta) (x- x_{\ell_2}) + y_{\ell_2} + \frac{2\Delta}{\cos(\theta)} \ge \tan(\beta) (x - x_{\ell_2}) + y_{\ell_2} +\frac{2\Delta}{\cos(\beta)} \quad &\forall \beta \in [\theta, \theta + \alpha]
        \end{aligned}
        \label{pr:inequality_ub1}
    \end{equation}

\item $x \ge x^I_u$:
\begin{equation}
        \begin{aligned}
            x&\ge x_{\ell_2} - 2\Delta \left( \frac{\sqrt{1+\tan^2(\theta+\alpha)} - \sqrt{1+\tan^2(\theta)}}{\tan(\theta+\alpha) - \tan(\theta)} \right)
            \\
            & \ge x_{\ell_2} - 2\Delta \left( \frac{\sqrt{1+\tan^2(\theta+\alpha)} - \sqrt{1+\tan^2(\beta)}}{\tan(\theta+\alpha) - \tan(\beta) } \right) & \forall \beta \in [\theta, \theta + \alpha]
             \\ \implies
            & \tan(\theta+\alpha) (x- x_{\ell_2}) + y_{\ell_2} + \frac{2\Delta}{\cos(\theta+\alpha)} \ge \tan(\beta) (x - x_{\ell_2}) + y_{\ell_2} +\frac{2\Delta}{\cos(\beta)} \quad &\forall \beta \in [\theta, \theta + \alpha]
        \end{aligned}
        \label{pr:inequality_ub2}
    \end{equation}
\end{enumerate}
Hence, the upper bound of $I(x)$ will be (see Figure~\ref{fig:INtersections.jpg})
\begin{equation}
\label{pr: upperbound_I}
\begin{cases}
     \tan(\theta) (x- x_{\ell_2}) + y_{\ell_2} + \frac{2\Delta}{\cos(\theta)} \text{ (line }\lambda_1) & x \le x^I_u
    \\
    \tan(\theta+\alpha) (x- x_{\ell_2}) + y_{\ell_2} + \frac{2\Delta}{\cos(\theta+\alpha)} \text{ (line }\lambda_3) & x \ge x^I_u
\end{cases}
\end{equation}

\begin{figure}[h!]
    \centering
    \includegraphics[width = 6cm]{ 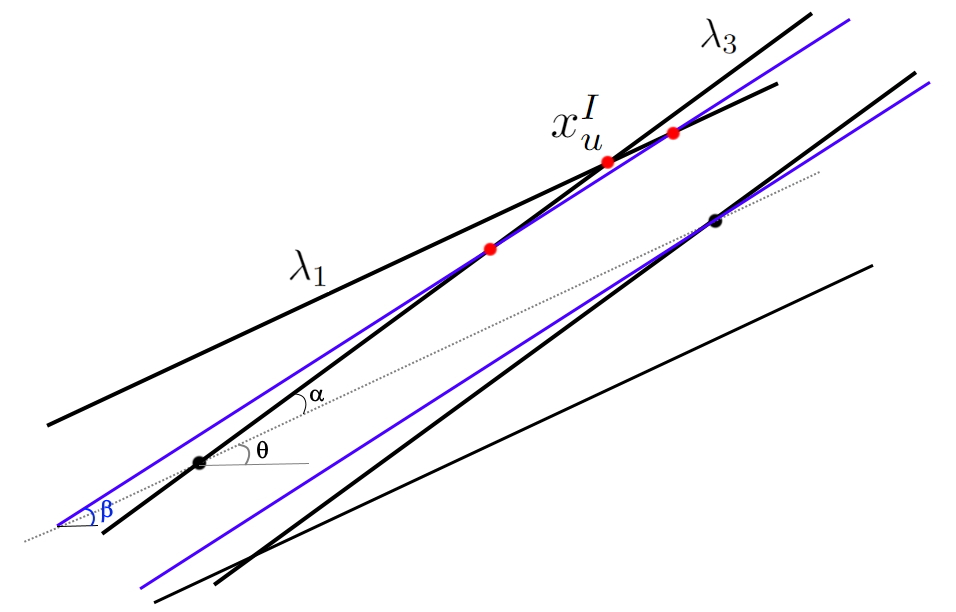}
    \caption{For $x \le x^I_u$, the boundary line $ \lambda_1$ and for $x \le x^I_u$, the boundary line $\lambda_3$ restrict any given line $\tan(\beta) (x - x_{\ell_2}) + y_{\ell_2} +\frac{2\Delta}{\cos(\beta)}, \forall \beta \in [\theta-\alpha, \theta+\alpha]$.}
\label{fig:INtersections.jpg}
\end{figure}

Likewise, upper bound of $J(x)$ can be determined by defining:
$$
x^J_u = x_{\ell_1} - 2\Delta \left( \frac{\sqrt{1+\tan^2(\theta)} - \sqrt{1+\tan^2(\theta-\alpha)}}{\tan(\theta) - \tan(\theta-\alpha)} \right),
$$
as the intersection of lines $\lambda_1$ and $\lambda_5$. Employing Lemma~\ref{lemma:increasing Function}, similarly it is concluded that the upper bound of $J(x)$ will be
\begin{equation}
\label{pr: upperbound_J}
\begin{cases}
     \tan(\theta) (x- x_{\ell_1}) + y_{\ell_1} + \frac{2\Delta}{\cos(\theta)} \text{ (line }\lambda_1) & x \ge x^J_u
    \\
    \tan(\theta-\alpha) (x- x_{\ell_1}) + y_{\ell_1} + \frac{2\Delta}{\cos(\theta-\alpha)} \text{ (line }\lambda_5) & x \le x^J_u
\end{cases}
\end{equation}

The lower bounds of $I(x)$ and $J(x)$ can be determined with a similar approach. Define
$$
x^I_l = x_{\ell_1} + 2\Delta \left( \frac{\sqrt{1+\tan^2(\theta+\alpha)} - \sqrt{1+\tan^2(\theta)}}{\tan(\theta+\alpha) - \tan(\theta)} \right),
$$
and 
$$
x^J_l = x_{\ell_2} + 2\Delta \left( \frac{\sqrt{1+\tan^2(\theta)} - \sqrt{1+\tan^2(\theta-\alpha)}}{\tan(\theta) - \tan(\theta-\alpha)} \right).
$$
Then, the lower bound of $I(x)$ will be 
\begin{equation}
\label{pr:lowerbound_I}
\begin{cases}
     \tan(\theta) (x- x_{\ell_1}) + y_{\ell_1} - \frac{2\Delta}{\cos(\theta)} \text{ (line }\lambda_2) & x \ge x^I_l
    \\
    \tan(\theta+\alpha) (x- x_{\ell_1}) + y_{\ell_1} - \frac{2\Delta}{\cos(\theta+\alpha)} \text{ (line }\lambda_4) & x \le x^I_l,
\end{cases}
\end{equation}
and for the lower bound of $J(x)$, we have
\begin{equation}
\label{pr:lowerbound_J}
\begin{cases}
     \tan(\theta) (x- x_{\ell_2}) + y_{\ell_2} - \frac{2\Delta}{\cos(\theta)} \text{ (line }\lambda_2) & x \le x^J_l
    \\
    \tan(\theta-\alpha) (x- x_{\ell_2}) + y_{\ell_2} - \frac{2\Delta}{\cos(\theta-\alpha)} \text{ (line }\lambda_6) & x \ge x^J_l.
\end{cases}
\end{equation}

In conclusion, for given location $(x_{b_3},y_{b_3})$ if it is located above the upper bounds in~\eqref{pr: upperbound_I} and~\eqref{pr: upperbound_J} which means
\begin{equation}
\label{pr:upperbounds}
    \begin{cases}
        y_{b_3} > \tan(\theta) (x_{b_3}- x_{\ell_1}) + y_{\ell_1} + \frac{2\Delta}{\cos(\theta)} (= \tan(\theta) (x_{b_3}- x_{\ell_2}) + y_{\ell_2} + \frac{2\Delta}{\cos(\theta)}),
        \\
        y_{b_3} > \tan(\theta+\alpha) (x_{b_3}- x_{\ell_2}) + y_{\ell_2} + \frac{2\Delta}{\cos(\theta+\alpha)} (= \tan(\theta+\alpha) (x_{b_3}- x_{\ell_1}) + y_{\ell_1}),
        \\
        y_{b_3} > \tan(\theta-\alpha) (x_{b_3}- x_{\ell_1}) + y_{\ell_1} + \frac{2\Delta}{\cos(\theta-\alpha)} (= \tan(\theta-\alpha) (x_{b_3}- x_{\ell_2}) + y_{\ell_2}),
    \end{cases}
\end{equation}
\textbf{OR} below the lower bounds in \eqref{pr:lowerbound_I} and~\eqref{pr:lowerbound_J} which implies
\begin{equation}
\label{pr:lowerbounds}
    \begin{cases}
        y_{b_3} < \tan(\theta) (x_{b_3}- x_{\ell_1}) + y_{\ell_1} - \frac{2\Delta}{\cos(\theta)} (= \tan(\theta) (x_{b_3}- x_{\ell_2}) + y_{\ell_2} - \frac{2\Delta}{\cos(\theta)}),
        \\
        y_{b_3} < \tan(\theta+\alpha) (x_{b_3}- x_{\ell_1}) + y_{\ell_1} - \frac{2\Delta}{\cos(\theta+\alpha)} (= \tan(\theta+\alpha) (x_{b_3}- x_{\ell_2}) + y_{\ell_2}),
        \\
        y_{b_3} < \tan(\theta-\alpha) (x_{b_3}- x_{\ell_2}) + y_{\ell_2} + \frac{2\Delta}{\cos(\theta-\alpha)} (= \tan(\theta-\alpha) (x_{b_3}- x_{\ell_1}) + y_{\ell_1}),
    \end{cases}
\end{equation}
then $y_{b_3} \notin I(x)\cup J(x)$. The set of infeasible locations for $\ell_1$ and $\ell_2$, $\Gamma_{\ell_1,\ell_2}$, includes locations in $L$ that meet~\eqref{pr:upperbounds} or~\eqref{pr:lowerbounds}, i.e.,
\begin{equation}
\begin{aligned}
    \Gamma_{\ell_1,\ell_2} &= \{b_i \in L\setminus\{\ell_1,\ell_2\}: y_{b_i} \notin I(x_{b_i})\cup J(x_{b_j}) \}
    \\
    &=\{b_i \in L\setminus\{\ell_1,\ell_2\}: (x_{b_i}, y_{b_i}) \in R''_{\ell_1,\ell_2} \}. 
\end{aligned}
\end{equation}
Therefore, locations in $\Gamma_{\ell_1,\ell_2}$ cannot be covered by a line that has a distance at most $\Delta$ from two locations $\ell_1$ and $\ell_2$, and so $z_{\ell_1}+z_{\ell_2}+z_{b_3}\le 1$ for $b_3 \in \Gamma_{\ell_1,\ell_2}$ is a valid cut.

\begin{remark}
The boundaries of the cases $\theta+\alpha \ge \pi/2$ or $\theta-\alpha \le -\pi/2$ can be easily obtained by rotating coordinates to get $\theta'$ such that $\theta'+\alpha < \pi/2$ and $\theta'-\alpha > -\pi/2$ the same as the conditions in the proof. 
\end{remark}\epf

\section{Details on the Feasibility Check when $E<\infty$}\label{ssec:feasibility}

Next we detail the feasibility check of the DBC for a given integer solution $z^h$ for the case where $E<\infty$ and that the solution of the SLA generates a line of length greater than $E$. Let $C\subseteq L$ be a subset of locations in $L$ and let $\lambda$ be a line in $\Rbb^2$. We say that $\lambda$ is a \emph{covering line} for $C$ if the distance from $\lambda$ to any $\ell\in C$ is at most $\Delta$. Similarly, we will refer to a finite segment $\xi$ that covers all locations in $C$ a \emph{covering segment} for $C$; the length of segment $\xi$ will be denoted by $u_\xi$. For any $a,b\in L$, let $\lambda_a$ and $\lambda_b$ the projections of $a$ and $b$ in $\lambda$, respectively.

\begin{lemma}\label{lemma:projection1}
Let $C\subseteq L$ be given and let $\lambda$ be a covering line for $C$. If $D(\lambda_a,\lambda_b)\le E$ for all $a,b\in C$, then there exists a covering segment $\xi$ of $C$ such that $u_\xi\le E$.
\end{lemma}
\noindent\emph{Proof:} Let $w,w'\in\argmax\{D(v,v')\colon v=\lambda_a,v'=\lambda_b, a,b\in C\}$. That is, $w$ and $w'$ are the projections of the points in $C$ on $\lambda$ that attain the maximum distance. If $\xi=[w,w']$, then clearly $\xi$ covers $C$ and $u_\xi\le E$.\epf

Let $C\subseteq L$ be given and let $\lambda$ be a covering line of $C$. We denote by $C_{\lambda,E}\subseteq C$ the set of building pairs in $C$ whose projections in $\lambda$ are at a distance of more than $E$, i.e., 
\begin{equation}
    C_{\lambda,E}=\bigl\{(a,b)\in C\times C\colon D(\lambda_a,\lambda_b)>E\bigr\}.
\end{equation}
For any $\ell\in L$ let $\delta_\ell=D(\ell,\lambda_\ell)$ and let $\xi_\ell=\sqrt{\Delta^2-\delta_\ell}$. For any $(a,b)\in C_{\lambda,E}$ let $\lambda_a'=\lambda_a+\xi_a(\lambda_b-\lambda_a)$ and $\lambda_b'=\lambda_b+\xi_b(\lambda_a-\lambda_b)$. Note that both $\lambda_a'$ and $\lambda_b'$ are points in $\lambda$.

\begin{lemma}\label{lemma:projection2}
Let $C\subseteq B$ be given and let $\lambda$ be a covering line for $C$. If $D(\lambda_a',\lambda_b')\le E$ for all $a,b\in C_{\lambda,E}$ then there exist a covering segment $\xi$ of $C$ such that $u_\xi\le E$.
\end{lemma}

\noindent\emph{Proof:} By construction, for any $(a,b)\in C_{\lambda,E}$, both $\lambda_a'$ and $\lambda_b'$ are elements of $\lambda$ and  $D(\lambda_a',a)=D(\lambda_b',b)=\Delta$, thus the segment $[\lambda_a',\lambda_b']$ covers both $a$ and $b$. Let $w,w'\in\argmax\{D(\lambda_a',\lambda_b')\colon a,b\in C_{\lambda,E}\}$, then by definition $D(\lambda_a',\lambda_b')\le D(w,w')$ for all $a,b\in C_{\lambda,E}$ and, moreover, $[w,w']$ covers $C$. Because  $D(w,w')\le E$ then $\xi=[w,w']$ is a covering segment of $C$ with $u_\xi\le E$. \epf

Step~\ref{algStep: check feasibility} in Algorithm~\ref{alg:Separation_Procedure} which checks the feasibility of solution $z^h$ in $\Uscr$, can be reframed as finding two endpoints $e_0$ and $e_1$ of some covering segment $\xi$ for $C^h$ with $u_\xi \le E$.
To accelerate the search for two endpoints in step \ref{algStep: check feasibility}, we implement SLA for $C^h$ along with two lemmas~\ref{lemma:projection1} and~\ref{lemma:projection2} as shown in Algorithm~\ref{feasibility_check_steps}.

\begin{algorithm}
    \KwData{Set $\Tilde{\Uscr}^h,C^h,E$}
    \KwResult{$feasible$}
    set $feasible = FALSE$\;
    run SLA for locations in $C^h$. Let $\lambda_{MAX}$ be the line intersecting maximal number of circles\label{step:stabbing_line} and $\xi_{MAX}$ be the shortest covering segment\;
    \If{maximum coverage with $\lambda_{MAX} = |C^h|$}{
    \If{$u_{\xi_{MAX}}\le E$}{\label{step:segment_stabLine}
    $feasible = TRUE$
    }
    \Else{
    \If{$\Tilde{\Uscr}^h \neq \emptyset$}{\label{step:feasiblity_set}
        $feasible = TRUE$
    }
    }
    }
    \caption{\label{alg:check}Feasibility check for solution $z^h$ if $E<\infty$}\label{feasibility_check_steps}
\end{algorithm}

Algorithm~\ref{feasibility_check_steps} executes SLA for the set of circles with centers in $C^h$ and radius $\Delta$ (step~\ref{step:stabbing_line}). If the maximal stabbing line $\lambda_{MAX}$ does not cover locations in $C^h$, then no segment produces the exiting solution which means $z^h$ is infeasible. Else, using the discussion in Lemmas~\ref{lemma:projection1} and~\ref{lemma:projection2}, we check whether the shortest covering segment of $\lambda_{MAX}$, denoted as $\xi_{MAX}$, has a length at most $E$ (step~\ref{step:segment_stabLine}). If so, we recognize $z^h$ feasible by $\xi_{MAX}$. If the length of the segment is greater than $E$, we find a feasible solution in $\Tilde{\Uscr}^h$. If set $\Tilde{\Uscr}^h$ is not empty, then there are feasible endpoints for a segment covering $C^h$ (step~\ref{step:feasiblity_set}). 

In step~\ref{step:feasiblity_set} we further enhance the tractability of $\Tilde{\Uscr}^h$ by restricting the search in ranges for feasible endpoints as follows. Let $\ell_i$ and $\ell_j$ be two extreme locations that have the maximum distance among all pairs of locations in $C^h$. Then, the endpoints of the line segment $e_0$ and $e_1$ must be within distance $\Delta$ from each $\ell_i$ and $\ell_j$. Particularly, points $e_0$ and $e_1$ should satisfy that (see Figure~\ref{fig:NewConstraintforUncertainty}):
\begin{subequations}\label{eq:linearcuts}
    \begin{align}
        x_{\ell_i} -\Delta \le e^x_0 \le x_{\ell_i} +\Delta
        \\
        y_{\ell_i} -\Delta \le e^y_0 \le y_{\ell_i} +\Delta
        \\
        x_{\ell_j} -\Delta \le e^x_1 \le x_{\ell_j} +\Delta
        \\
        y_{\ell_j} -\Delta \le e^y_1 \le y_{\ell_j} +\Delta.
    \end{align}
\end{subequations}
The enhancement, therefore, adds the inequalities~\eqref{eq:linearcuts} in the definition of $\Tilde{\Uscr}^h$. The enhanced definition is solved significantly faster.

\begin{figure}[h]
    \centering
    \includegraphics[width=6cm]{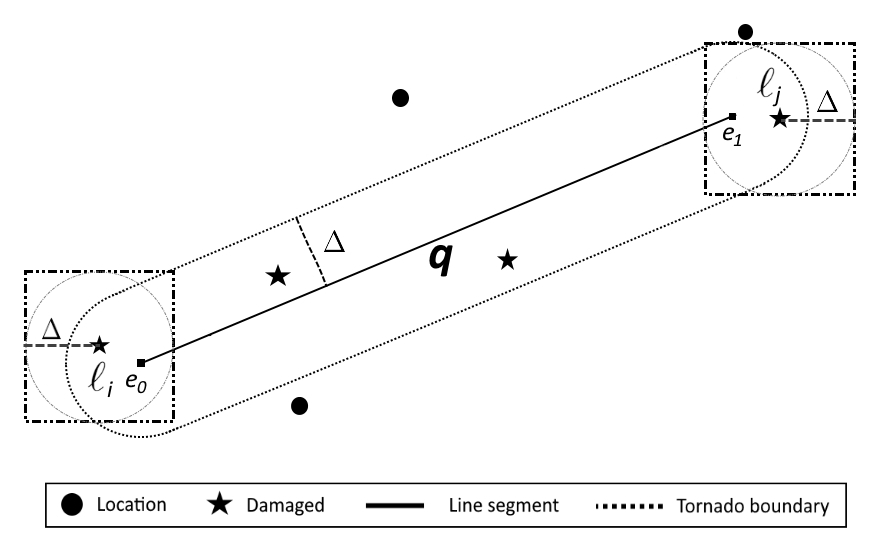}
    \caption{Two endpoints $e_0$ and $e_1$ are within circles centered at locations $\ell_i$ and $\ell_j$ with radius $\Delta$. These two circles are inside squares which are added to restrict the search ranges for endpoints by linear constraints.}
    \label{fig:NewConstraintforUncertainty}
\end{figure}

\section{Performance of Subproblem Methods}\label{sec:performance-compare}
We compare three different methods to address the one-level subproblem $\Phi(f)$ in Algorithm~\ref{alg:algorithm_seq}:
\begin{enumerate}
    \item ORG: solves the subproblem with \textit{original} set $\Uscr$ (instead of using \eqref{eq:infeasible_comb_indices}) by DBC and only updates $\mathcal{R}(f)$ on the fly,
    \item AVC: \textit{adds valid cuts} for infeasible combinations in $\mathcal{C}^0$ to $\Uscr$ to enhance ORG performance while updates $\mathcal{R}(f)$ on the fly,
    \item DEC: replace $\Uscr$ by conflict constraints~\eqref{eq:infeasible_comb_indices} and update $\mathcal{C}$ and $\mathcal{R}(f)$ on the fly as described in Section~\ref{subsec:DBCoverview}.
\end{enumerate}
Table~\ref{Tab1:comparison} compares the three approaches in terms of optimality gap, solution time, number of iterations, subproblem run time, and its callback function run time. We generate a testbed of 5 sample problems in which 10 blocks in Joplin are randomly selected as locations. Column 2 shows the method that is implemented to solve the subproblem. Columns 3, 4, and 5 present the best objective value (expected population dislocation), the best lower bound, and the optimality gap found within the time limit of one hour, respectively. Columns 6 and 7 show the CPU run time of C\&CG Algorithm~\ref{alg:algorithm_seq} in seconds and its number of iterations, respectively. Column 8 compares the CPU run time of the DBC method in seconds for each selection of subproblems. Column 9 also presents the total amount of time that the DBC spends in the callback function to generate cuts on the fly with Algorithm~\ref{alg:Separation_Procedure}; note that the callback run time for ORG and AVC methods does not include the time of checking feasibility in $\Uscr$ as we directly use the original feasible solution.

\begin{table}[!ht]
    \centering
    \begin{adjustbox}{max width=\textwidth}
    \begin{tabular}{c|c|ccc|cc|cc}
    \hline
        Sample \# & Subproblem method & Best objective & Best bound & Gap & CCG run time (sec.) & CCG iteration & DBC run time (sec.) & Callbacks run time (sec.) \\ \hline
        \multirow{3}{*}{1} & ORG & 365 & - & - & \textit{3600} & 1 & 3600 & 42 \\
        ~ & AVC & 185 & 185 & 0\% & 6 & 3 & 6 & 0 \\ 
        ~ & DEC & 185 & 185 & 0\% & 1 & 3 & 1 & 1 \\ \hline
        \multirow{3}{*}{2} & ORG & 161 & 153 & 5\% & \textit{4012} & 2 & 4012 & 1 \\
        ~ & AVC & 153 & 153 & 0\% & 26 & 3 & 26 & 1 \\
        ~ & DEC & 153 & 153 & 0\% & 1 & 3 & 1 & 0 \\ \hline
        \multirow{3}{*}{3} & ORG & 178 & 178 & 0\% & \textit{4484} & 3 & 4484 & 1 \\
        ~ & AVC & 178 & 178 & 0\% & 9 & 3 & 9 & 0 \\ 
        ~ & DEC & 178 & 178 & 0\% & 2 & 3 & 1 & 1 \\ \hline
        \multirow{3}{*}{4} & ORG & 56 & - & - & \textit{3600} & 1 & 3600 & 0 \\ 
        ~ & AVC & 75 & 75 & 0\% & 3549 & 3 & 3549 & 1 \\
        ~ & DEC & 75 & 75 & 0\% & 1 & 3 & 1 & 1 \\ \hline
        \multirow{3}{*}{5} & ORG & 103 & - & - & \textit{3600} & 1 & 3600 & 0 \\ 
        ~ & AVC & 65 & 65 & 0\% & 36 & 3 & 36 & 1 \\ 
        ~ & DEC & 65 & 65 & 0\% & 1 & 3 & 1 & 0 \\ \hline
    \end{tabular}
    \end{adjustbox}
       \caption{Comparing C\&CG performance for 3 different methods to solve the subproblem}
       \label{Tab1:comparison}
\end{table}

The results in Table~\ref{Tab1:comparison} show that the ORG method is not able to solve any sample with only 10 locations within the one-hour time limit. Also, based on the results, solving the master problem and the second stage recovery problem in the callback function is very fast and almost all the run time in ORG is spent on solving the MINLP subproblem~\eqref{eq:tornado_bilevel_problem}. For samples 1, 4, and 5 the subproblem was not solved to optimality in one hour. For some samples, the ORG method exceeds the time limit while solving the MINLP subproblem. In these cases, we wait until the subproblem finishes its execution before terminating the algorithm, which results in some computational times reported being larger than the one-hour time limit (\textit{italic values} in the table).

Adding infeasible pairs and triples cuts in AVC significantly reduces the computation times by as much as two orders of magnitude, resulting in all the samples being solved in less than one hour. However, AVC can be potentially time-consuming for problems such as sample 4 in the experiment. The DEC method greatly outperforms the other two methods as it solved all samples within a few seconds. The subproblem runtime for DEC shows how fast the subproblem can lead to finding optimal tornado paths. The callbacks' runtime also emphasizes that the feasibility check with Algorithm~\ref{feasibility_check_steps} in Appendix~\ref{ssec:feasibility} which is equipped by the polynomial stabbing line algorithm is a quick task.

\section{\textcolor{black}{Additional Experiments for 30 Days of Recovery}}\label{sec:30daysRecoveryResults}
\textcolor{black}{We conduct more experiments for locations in Joplin to study the population dislocation after 30 days of recovery. The results are shown in Table~\ref{Tab:JoplinRobust-30days}.}

\begin{table}[!ht]
    \centering
    \captionsetup{justification=centering}
    \caption{Solving the two-stage robust optimization model with different parameters of tornado length and available budget,  after 30 days of recovery for 100 locations group in Joplin, MO.}
    \begin{adjustbox}{max width=\textwidth}
    \begin{tabular}{c|c|c|cc|ccc|cc}
    \hline
        Length & Budget & Population  & CCG  & CCG & DBC  & Callback  & $\Tilde{\mathcal{U}}^h$ feasibility & Retrofitting cost & Recovery cost\\ 
        & (\$) & dislocation & run time (sec.) & iteration  & run time (sec.) & run time (sec.) & run time (sec.) & /budget & /budget  \\
        \hline
        \multirow{3}{*}{5} & 0 M & 16465 & 153 & 2 & 153 & 5 & 2 & - & - \\
        ~ & 15 M & 14847 & 194 & 2 & 193 & 5 & 2 & 0.6 & 0.4 \\
        ~ & 30 M & 14688 & 162 & 2 & 161 & 5 & 1 & 0.3 & 0.7 \\ \hline
        \multirow{3}{*}{$\infty$} & 0 M & 17449 & 361 & 2 & 360 & 3 & 0 & - & - \\
        ~ & 15 M & 15742 & 527 & 2 & 527 & 5 & 0 & 0.6 & 0.4 \\
        ~ & 30 M & 15582 & 500 & 2 & 500 & 3 & 0 & 0.3 & 0.7 \\ \hline
    \end{tabular}
    \end{adjustbox}
    
       \label{Tab:JoplinRobust-30days}
\end{table}

\textcolor{black}{
Compared with the results in Table~\ref{Tab2:JoplinRobust}, we observe more population dislocation for each set of parameters because of the shorter period to recover the damaged locations.
Despite this difference in the magnitude of population dislocation, it is important to note that the information, clarifications, and examinations between both Table~\ref{Tab2:JoplinRobust} and~\ref{Tab:JoplinRobust-30days} remain consistent. Therefore, Table~\ref{Tab:JoplinRobust-30days} provides numerical experiments with a focus on the impact of shorter recovery times, complementing the information presented in Table~\ref{Tab2:JoplinRobust} in the main body of the paper.}

\end{document}